\documentstyle[12pt,psfig]{article}
\input{psfig}
\newtheorem{theorem}{Theorem}[section]
\newtheorem{lemma}[theorem]{Lemma}
\newtheorem{corollary}[theorem]{Corollary}

\newtheorem{problem}[theorem]{Problem}

\newtheorem{definition}[theorem]{Definition}
\newtheorem{example}[theorem]{Example}
\newtheorem{examples}[theorem]{Examples}
\newtheorem{question}[theorem]{Question}
\newtheorem{conjecture}[theorem]{Conjecture}

\begin{document}
\begin{center} Three talks in Cuautitlan under the general title \end{center}
\centerline{ {\bf Algebraic Topology Based on Knots }}
\centerline{ Topolog\'ia algebraica basada sobre nudos}
      \centerline{ March 12,13, 2001}
\centerline{ J\'ozef H.Przytycki }

\section{Talk 1: \ Open problems in knot theory that everyone can try 
to solve.}

Knot theory is more than two hundred years old; the first scientists 
who considered knots as mathematical objects were A.Vandermonde (1771) and
C.F.Gauss (1794). However, despite the impressive grow of the theory,
there are simply formulated but fundamental questions, to which we do not
know answers. I will discuss today several such open problems, describing
in detail the 20 year old Montesinos-Nakanishi
conjecture. Our problems lead to sophisticated mathematical structures
(I will describe some of them in tomorrows talks), but today's description
will be absolutely elementary.
Links are circles embedded in our space, $R^3$, up to topological deformation,
that is two links are equivalent if one can be deformed to the other 
in space without cutting and pasting. We represent links on a plane using
their diagrams (we follow the terminology of Lou Kauffman's talk). 

First we should introduce the concept of an $n$ move on 
a link. One should stress that the move is not topological, so that it can
change the type of the link we deal with.
\begin{definition}\label{1.1}
An $n$-move on a link is a local change of the link illustrated in 
Figure 1.1.
\end{definition}
\centerline{\psfig{figure=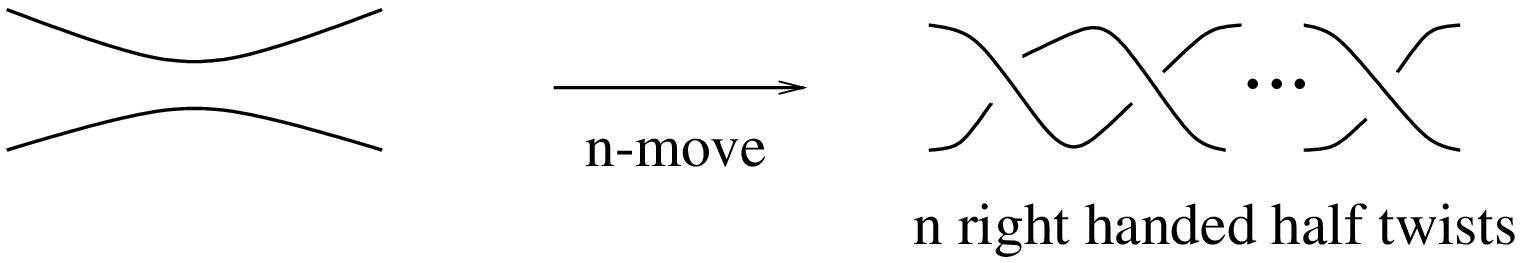,height=1.6cm}}
\centerline{Figure 1.1; \ $n$-move}
In our convention the part of the link outside of the disk in which the move
takes part, is unchanged. For example 
\psfig{figure=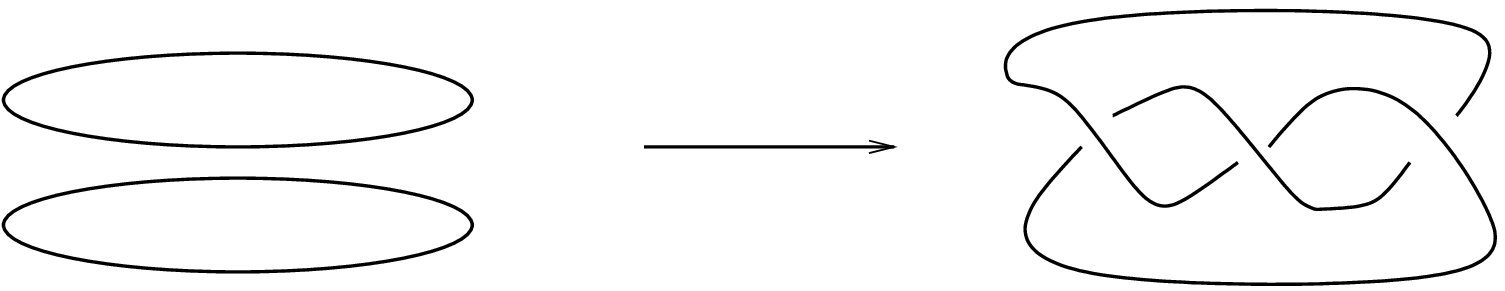,height=0.4cm}
illustrates a 3-move.
\begin{definition}\label{1.2}
We say that two links, $L_1, L_2$ are $n$-move equivalent if one can go 
from one to the other by a finite number of $n$-moves and their inverses
($-n$ moves). 
\end{definition}

If we work with diagrams of links then the topology of links is
reflected by Reidemeister moves, that is two diagrams represent
the same link in space if and only if one can go from one to the other
by Reidemeister moves:\ \\
\ \\
\centerline{\psfig{figure=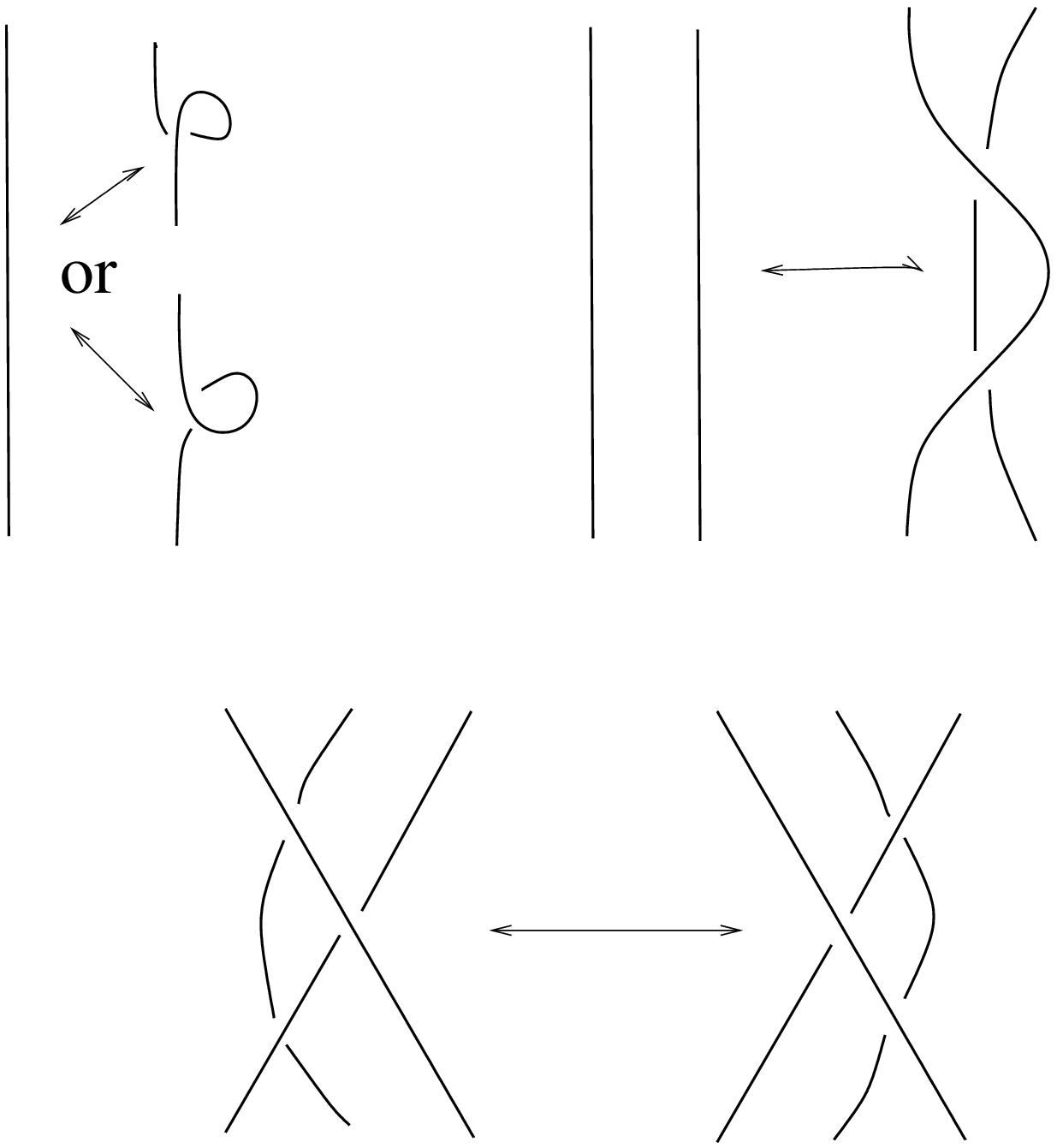,height=5.8cm}}
\centerline{Figure 1.2; \ Reidemeister moves}
Thus we say that two diagrams, $D_1$ and $D_2$, are $n$-move equivalent
if one can obtain one from the other by $n$-moves, their inverses and
Reidemeister moves. To illustrate this let us notice that the move
\psfig{figure=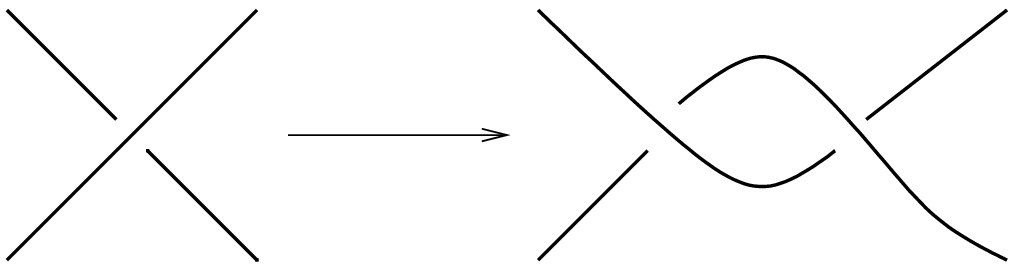,height=0.4cm}    
is a consequence of a 3-move and a second Reidemeister move (Fig.1.3).
\\
\ \\
\centerline{\psfig{figure=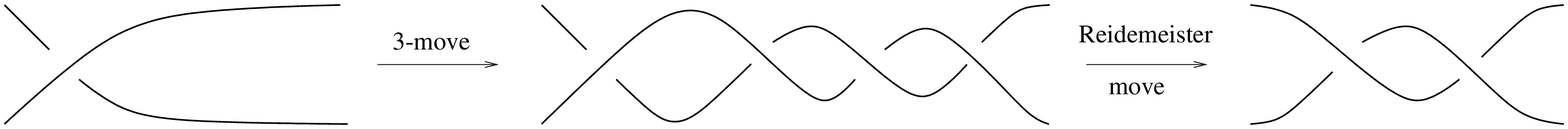,height=1.3cm}}
\centerline{Figure 1.3}

\begin{conjecture} [Montesinos-Nakanishi]
\label{Conjecture 1.3}\ \\
 Every link is $3$-move equivalent to a trivial link.
\end{conjecture}
Yasutaka Nakanishi first considered the conjecture in 1981. 
Jos\'e Montesinos analyzed 3-moves before, in connection with 
3-fold dihedral branch coverings, and asked a related but 
different question.

\begin{examples}\label{1.4}
\begin{enumerate}
\item [(i)] Trefoil knots are 3-move equivalent to the trivial
link of two components:\\
\ \\ 
\centerline{\psfig{figure=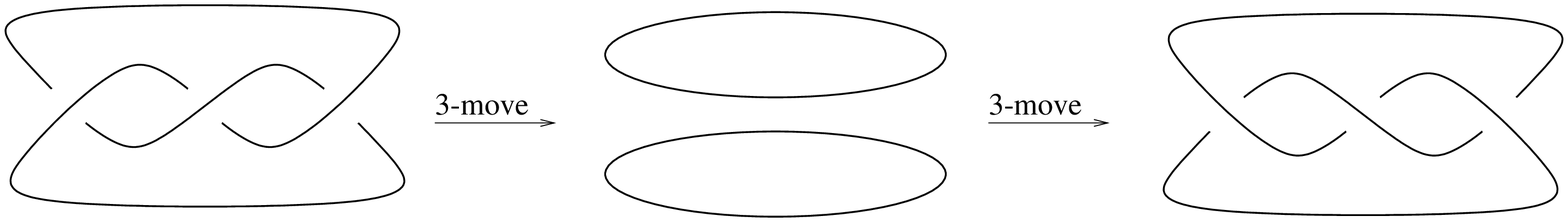,height=1.7cm}}
\centerline{Fig. 1.4} 
\item [(ii)] The figure eight knot is 3-move equivalent to the trivial
knot:\\
\ \\ 
\centerline{\psfig{figure=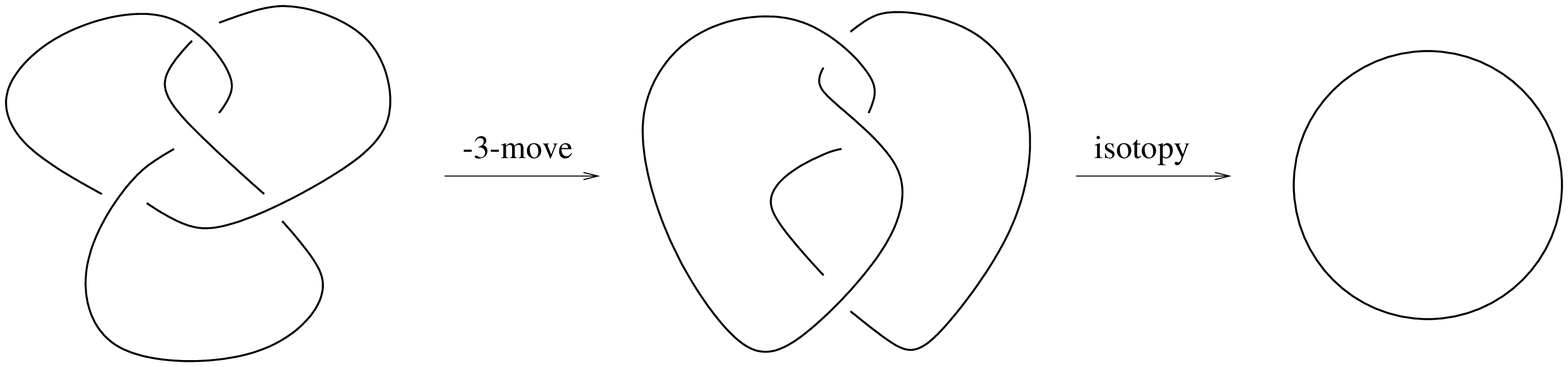,height=1.8cm}}
\centerline{Fig. 1.5}
\end{enumerate}
\end{examples}

We will show later, in this talk, that different trivial links
are not 3-move equivalent, however in order to achieve this
conclusion we need an invariant of links preserved by 3-moves and
different for different trivial links. Such invariant is the
Fox 3-coloring. We will introduce it later (today in the simplest
form and in the second lecture in a more general context of Fox
$n$-colorings and Alexander-Burau-Fox colorings). 
Now let us present some other related conjectures.

\begin{conjecture}\label{Conjecture 1.5}\ \\
 Any 2-tangle is 3-move equivalent, to one of the four 2-tangles
of Figure 1.6. We allow additionally trivial components in the tangles
of Fig.1.6.
\end{conjecture}
\ \\
\centerline{\psfig{figure=S+S-S0Sinf.eps}}
\begin{center}
Fig. 1.6
\end{center}

Montesinos-Nakanishi conjecture follows from Conjecture
1.5. More generally if Conjecture 1.5 holds for some class of 2-tangles,
then Conjecture 1.3 holds for a link obtained by closing any 
tangle from the class, without introducing any new
crossing. The simplest interesting tangles for which
Conjecture 1.5 holds  are algebraic tangles in the sense of Conway
(I will call them 2-algebraic tangles and present their
generalization in the next talk). For 2-algebraic tangles
Conjecture 1.5 holds by an induction and I will leave it
as a pleasent exercise for you. The necessary definition is given below:

\begin{definition}[\cite{Co,B-S}]\label{1.6}\ \\
2-algebraic tangles are the smallest family of 2-tangles which satisfies:\\
(0) Any 2-tangle with 0 or 1 crossing is 2-algebraic.\\
(1) If $A$ and $B$ are 2-algebraic tangles
    then $r^i(A)*r^j(B)$ is 2-algebraic; $r$ denotes here
the rotation of a tangle by $90^o$ angle along the $z$-axis, and * denotes
the (horizontal) composition of tangles.\\
A link is 2-algebraic if it is obtained from a 2-algebraic tangle
by closing its ends without introducing any new crossings.
\end{definition}

The Montesinos-Nakanishi conjecture was proven for many special
families of links by my students Q.Chen and T.Tsukamoto \cite{Che,Tsu,P-Ts}.
In particular Chen proved that the conjecture holds for all
five braid links except, possibly one family, represented by
the square of the center of the 5-braid group, 
$(\sigma_1\sigma_2\sigma_3\sigma_4)^{10}$. Chen found 5-braid link
3-move equivalent to it with 20 crossings. It is now the smallest known
possible counterexample to Montesinos-Nakanishi conjecture, Fig.1.7.  

\ \\
\centerline{\psfig{figure=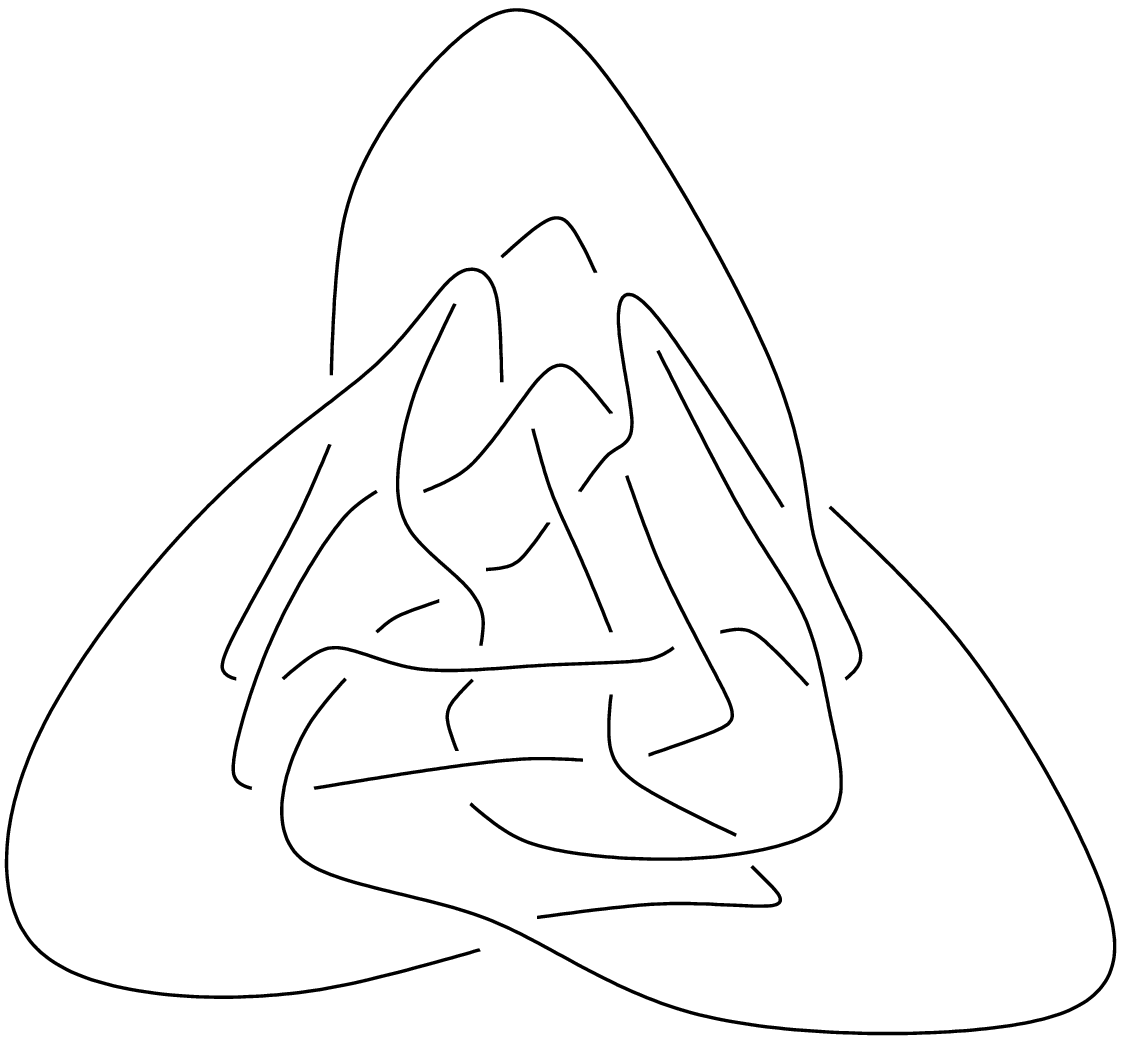,height=5.6cm}}
\begin{center}
Fig. 1.7
\end{center}
Previously (1994), Nakanishi suggested the 2-parallel of the
Borromean link (with 24 crossings) as a possible counterexample 
(Fig.1.8).

\ \\
\centerline{\psfig{figure=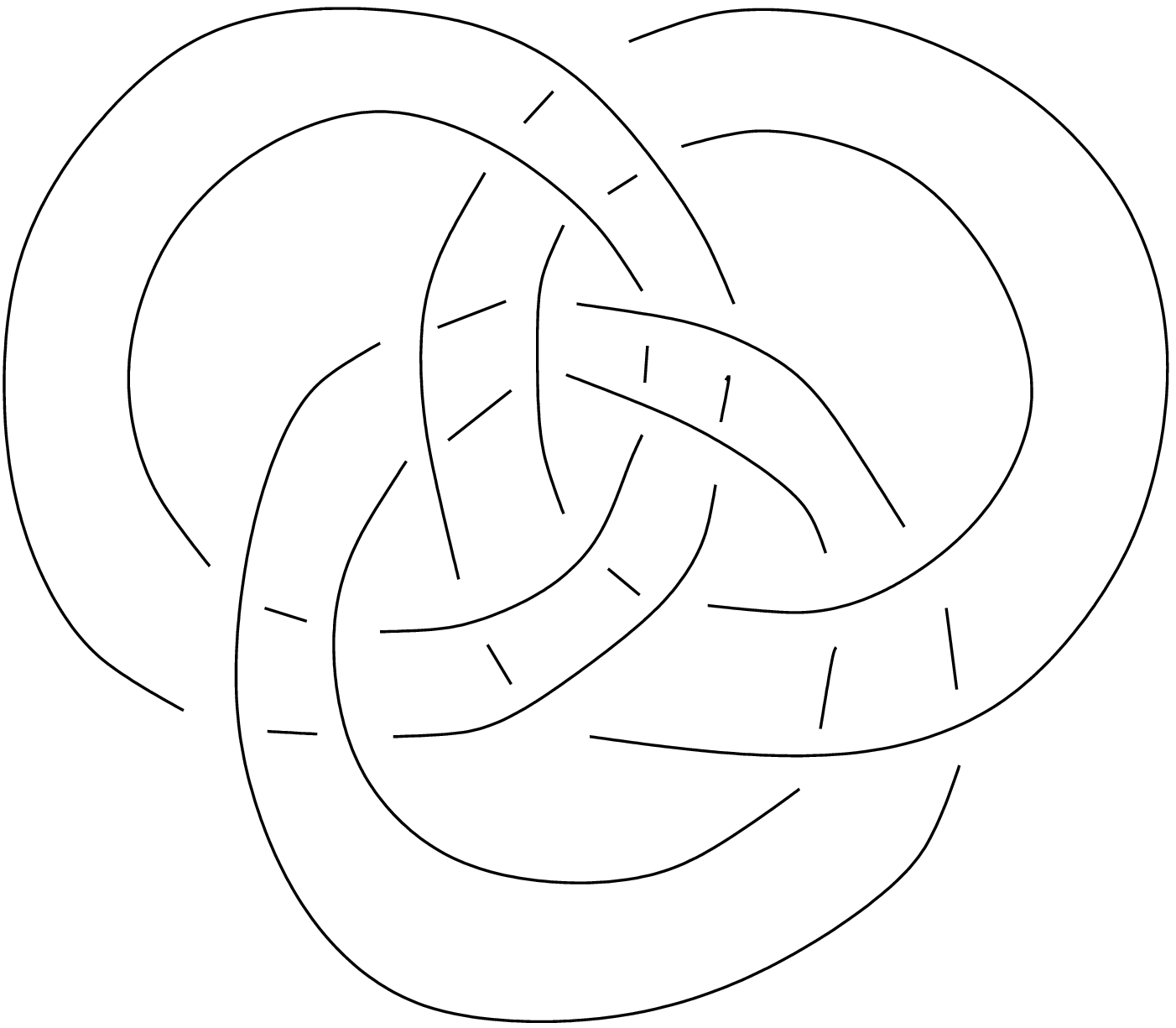,height=4.5cm}}
\begin{center}
Fig. 1.8
\end{center}

We will go back tomorrow to theories motivated by 3-moves,
now we will outline conjectures using other elementary moves.

\begin{conjecture}[Nakanishi, 1979]\label{1.7}\ \\
Every knot is $4$-move equivalent to the trivial knot.
\end{conjecture} 
\begin{examples}\label{1.8}
Reduction of the trefoil and the figure eight knot is 
illustrated in Figure 1.9.
\end{examples}
\centerline{\psfig{figure=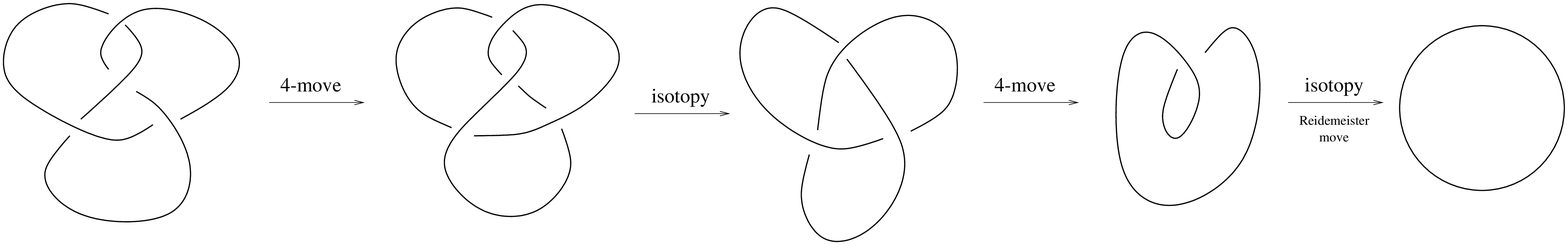,height=2.6cm}}
\centerline{Figure 1.9}

It is not true that every link is $5$-move equivalent to a trivial link.
One can show, using the Jones polynomial, that the figure eight knot is
not 5-move equivalent to any trivial link. One can however introduce
a more delicate move, called $(2,2)$-move 
(\psfig{figure=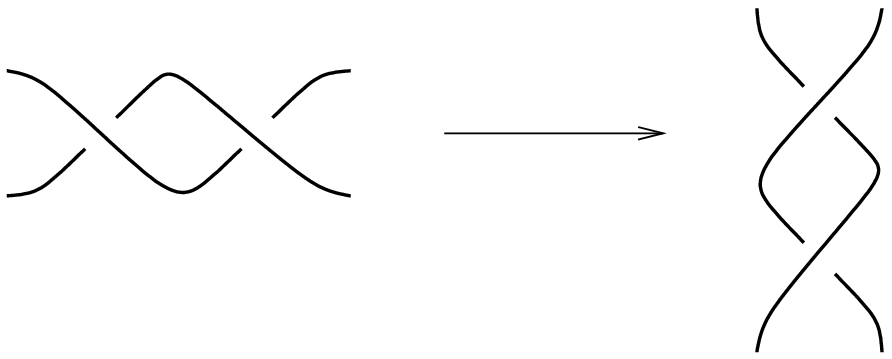,height=0.6cm})
 such that a $5$- move is a combination of a $(2,2)$-move 
and its mirror image $(-2,-2)$-move
(\psfig{figure=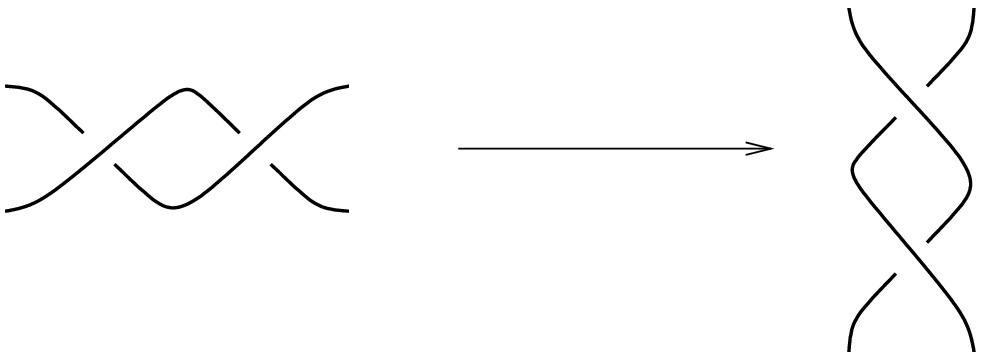,height=0.6cm}), as illustrated in 
Figure 1.10 \cite{H-U,P-3}.\\
\ \\

\centerline{\psfig{figure=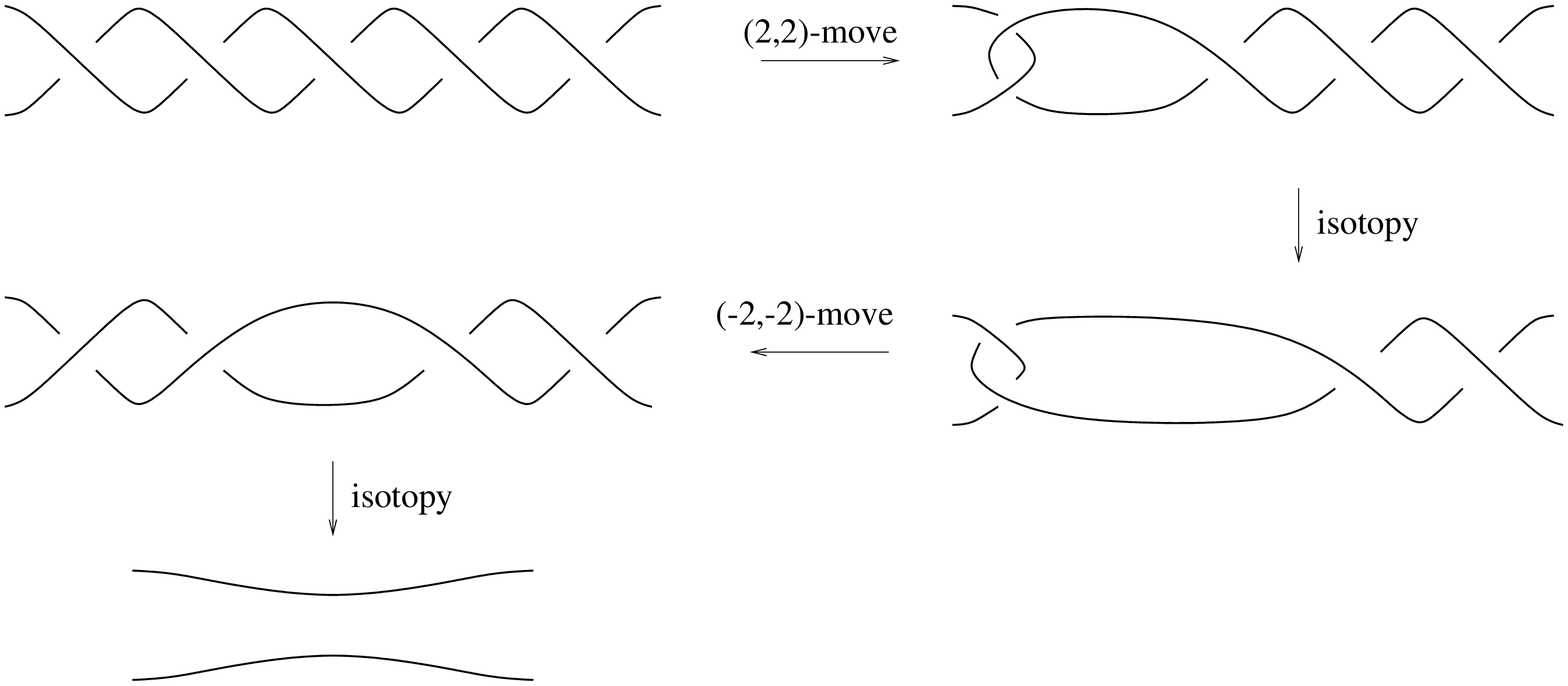,height=6.2cm}}
\centerline{Figure 1.10}

\begin{conjecture}[Harikae, Nakanishi, 1992]\label{1.9}
Every link is $(2,2)$-move equivalent to a trivial link.
\end{conjecture}
As in the case of 3-moves, an elementary induction shows that
the conjecture holds for 2-algebraic links 
(algebraic in the Conway's sense). It is also known for
all links up to $8$ crossings. The key element of the argument 
is the observation (going back to Conway \cite{Co}) that any link 
up to 8 crossings (different than $8_{18}$) is 2-algebraic.
The reduction of the $8_{18}$ knot to
a trivial link of two components by my students, Jarek Buczy\'nski and 
Mike Veve, is 
illustrated in Figure 1.11.
\ \\
\centerline{\psfig{figure=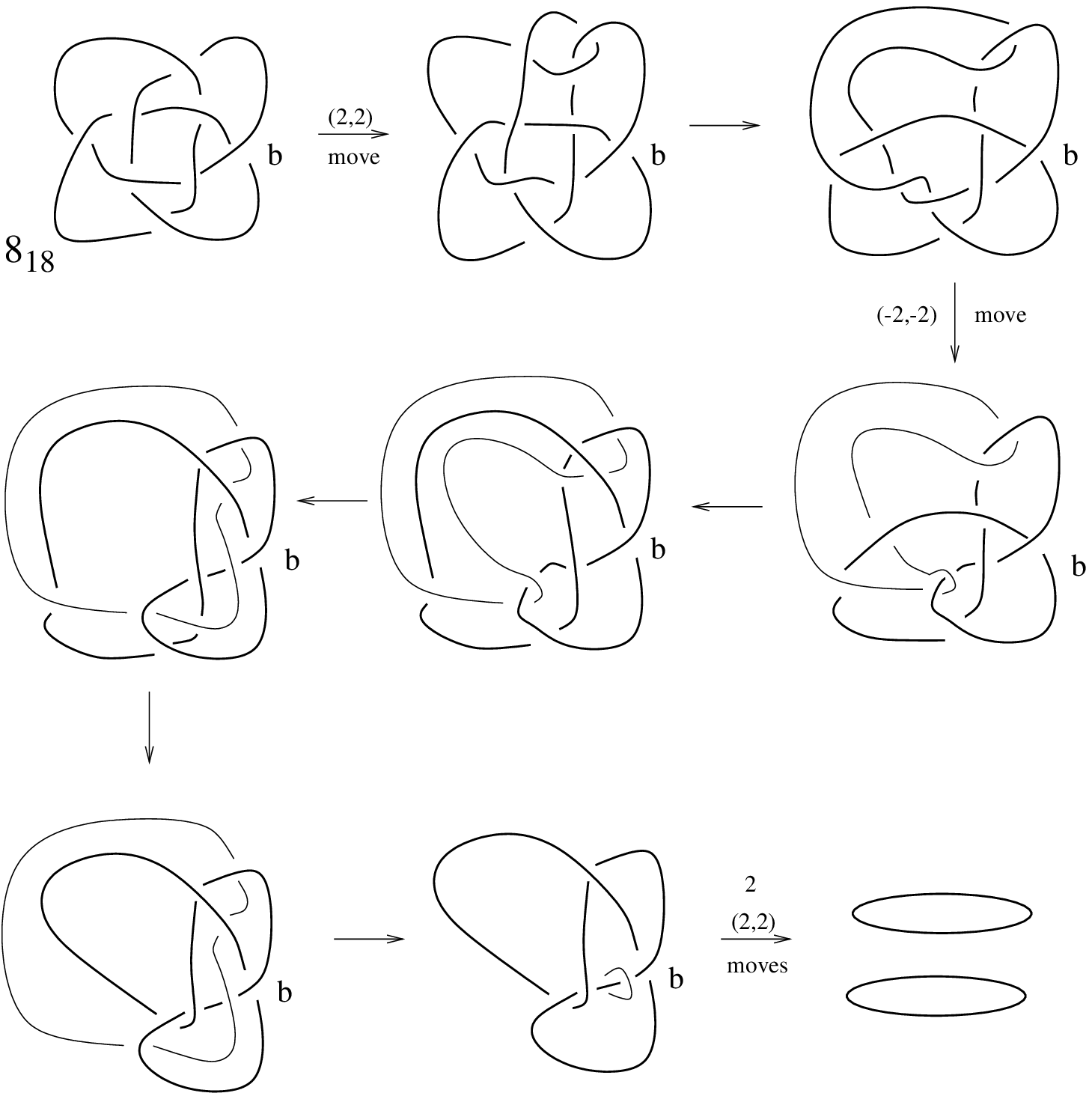,height=10.2cm}}

\centerline{Figure 1.11;\ Reduction of the $8_{18}$ knot}

The smallest knot, not reduced yet is the $9_{49}$ knot, Figure 1.12.
Possibly you can reduce it!
\ \\

\centerline{\psfig{figure=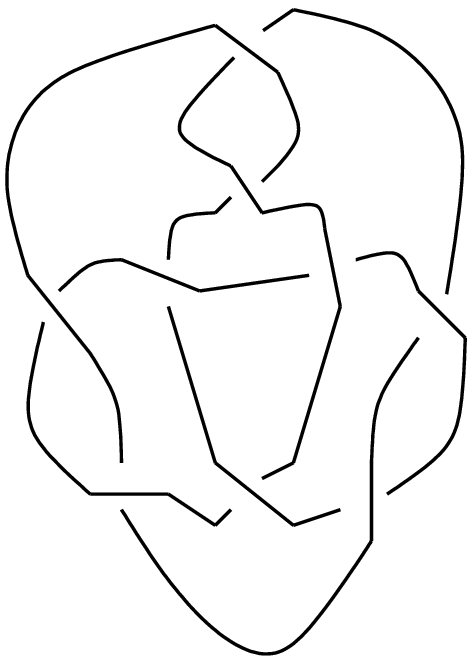,height=4.1cm}}
\centerline{Figure 1.12}

With the next open question, I am much less 
convinced that the answer is positive,
so I will not call it a conjecture.
First let define a $(p,q)$ move as a local modification of a link 
as shown in Figure 1.13.
We say that two links, $L_1,L_2$ are $(p,q)$ equivalent if one can go
from one to the other by a finite number of $(p,q)$,$(q,p)$,$(-p,-q)$ 
and $(-q,-p)$-moves.\\
\ \\
\centerline{\psfig{figure=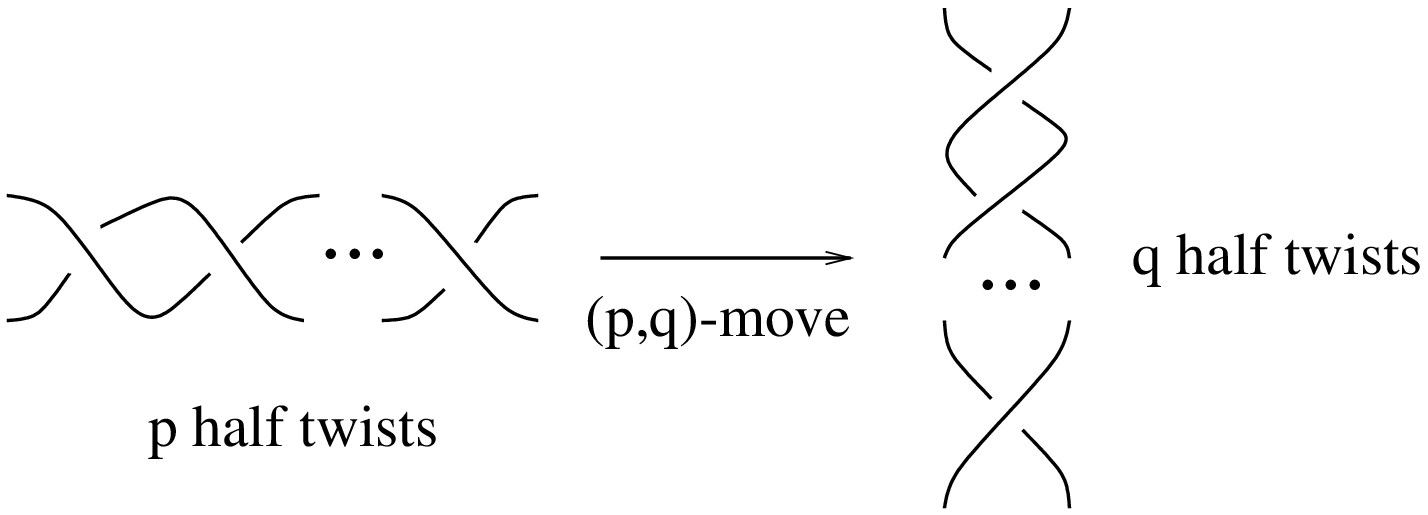,height=3.6cm}}
\centerline{Figure 1.13}

\begin{problem}[\cite{Kir};Problem 1.59(7), 1995]\label{1.10}
Is it true that any link is $(2,3)$ move equivalent to a trivial link?
\end{problem}

\begin{example}\label{1.11} Reduction of the trefoil and the figure eight knots is illustrated
in Fig. 1.14. Reduction of the Borromean rings is performed in Fig. 1.15.
\end{example}
\centerline{\psfig{figure=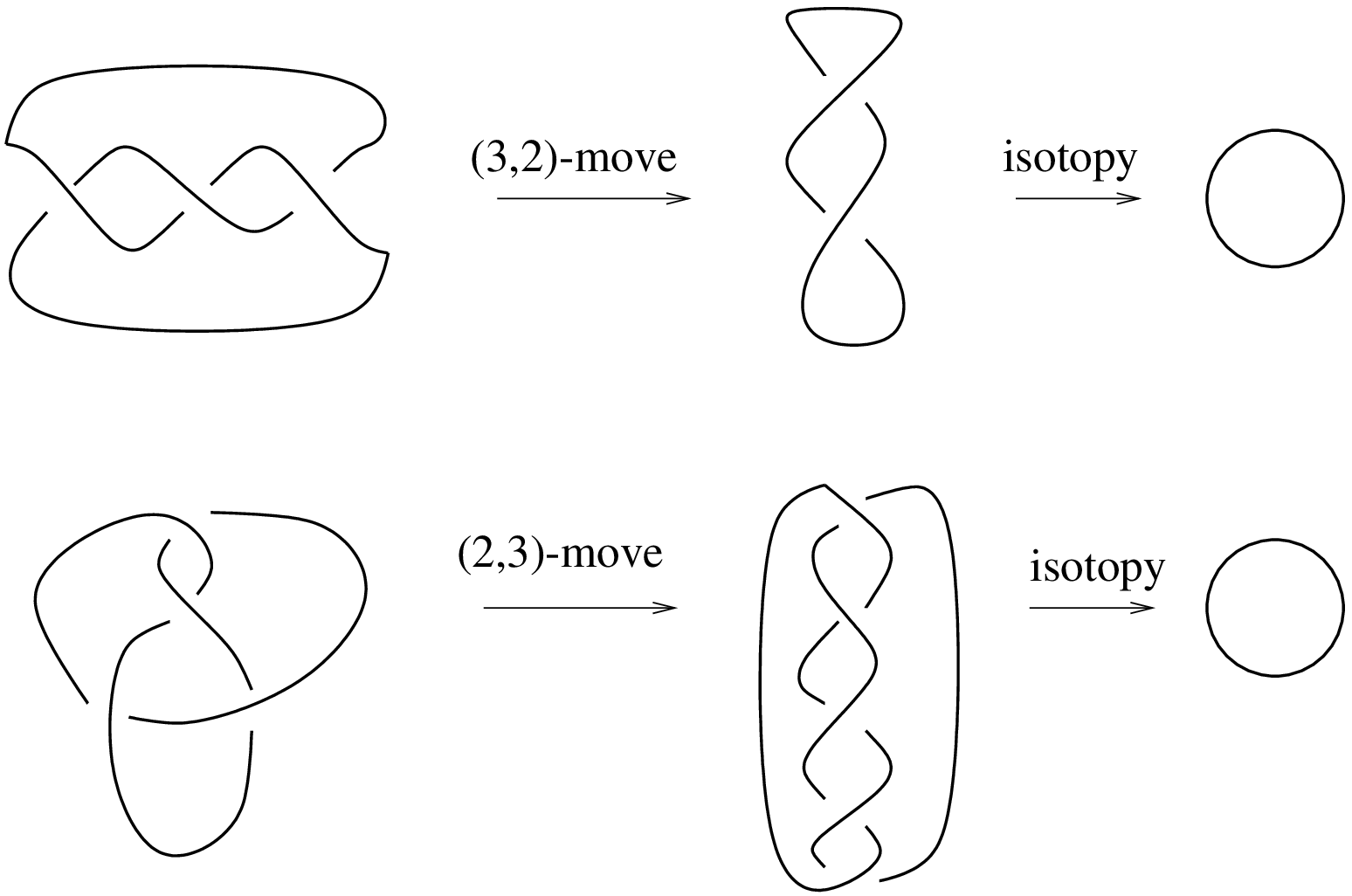,height=6.1cm}}
\centerline{Figure 1.14}

\centerline{\psfig{figure=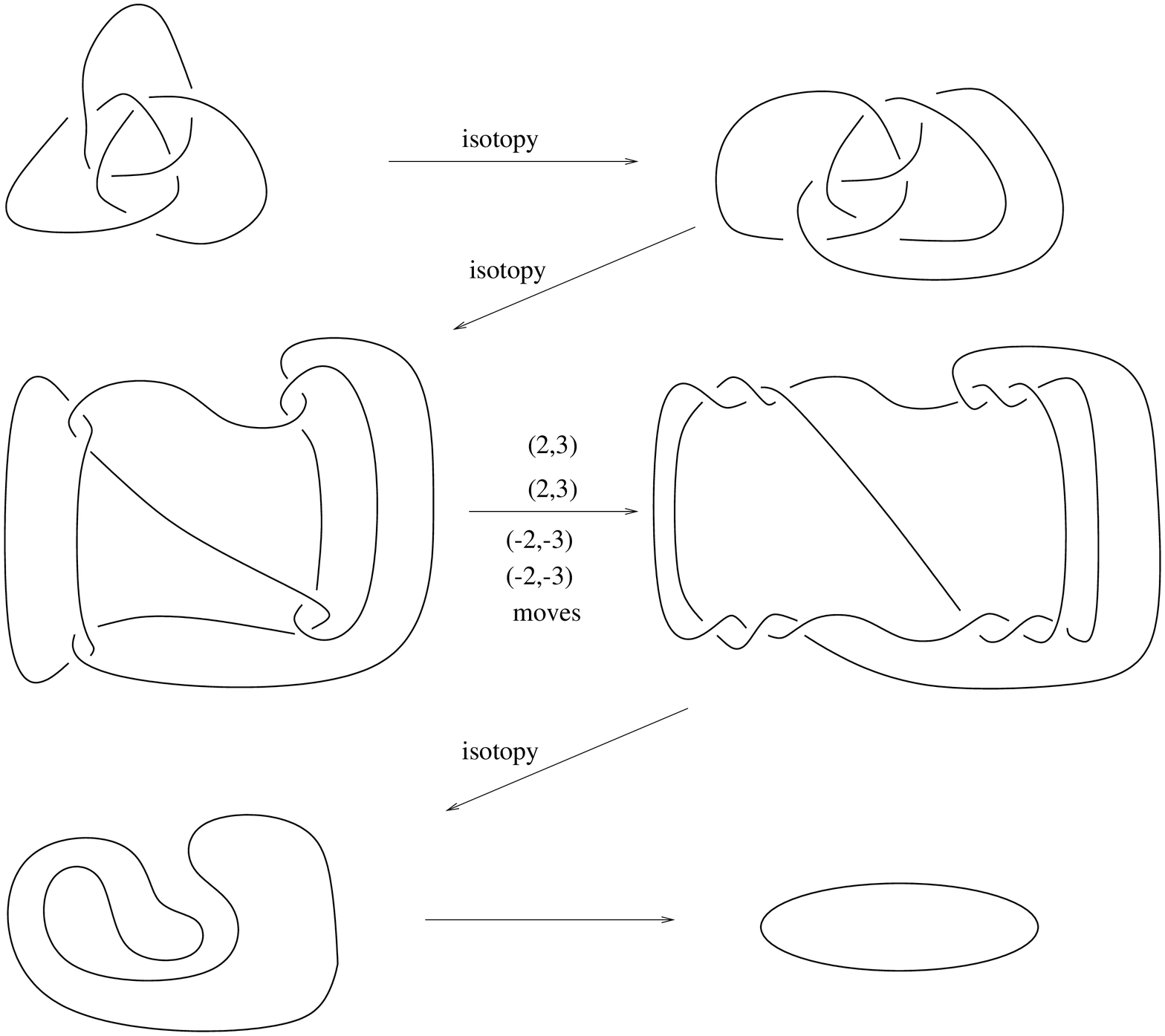,height=9.1cm}}
\centerline{Figure 1.15}

Generally, rather simple inductive argument shows that $2$-algebraic 
links are $(2,3)$-move equivalent to trivial links. 
Figure 1.16
illustrates why the Borromean rings are 2-algebraic. By properly
filling black dots one can also show that all links up to 8 crossings,
but $8_{18}$, are 2-algebraic. Thus,
as in the case of $(2,2)$-equivalence, the only link, up to $8$ crossings,
which should be still checked is the $8_{18}$ knot\footnote{To prove
that that the knot $8_{18}$ is not 2-algebraic
one would have to consider the 2-fold cover of $S^3$ with this knot
as a branching set and show that it not Waldhausen graph manifold.
In fact it is a hypebolic manifold so cannot be a graph manifold.}.
 Nobody really tried this
seriously, so maybe somebody in the audience will try this puzzle.\\
\ \\
\centerline{\psfig{figure=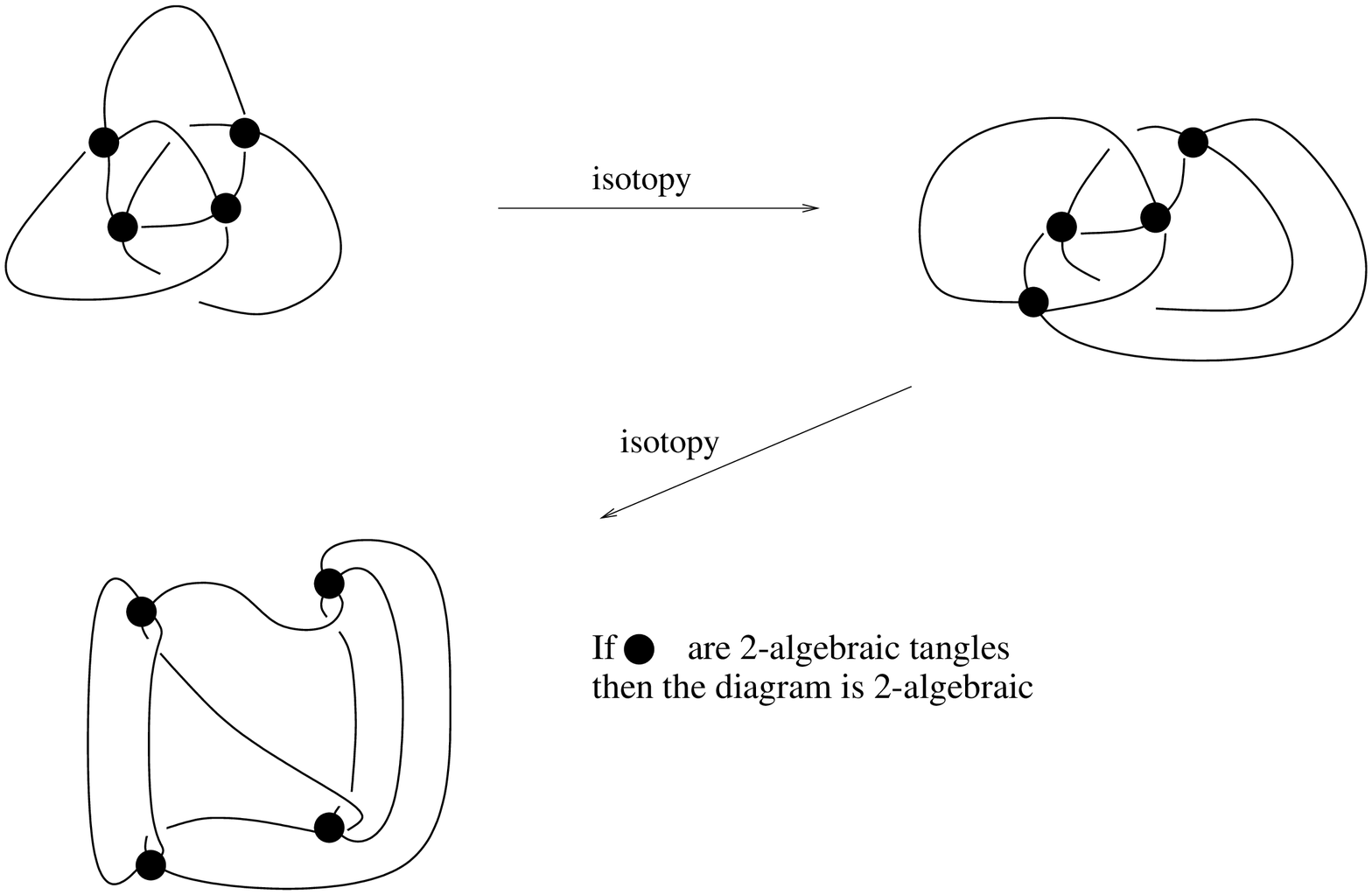,height=6.1cm}}
\centerline{Figure 1.16}
\ \\
{\bf Fox colorings.}\\
The 3-coloring which we will use to show that different 
trivial links are not 3-move equivalent,
was introduced by Ralph H. Fox in about 1956 when 
explaining knot theory to undergraduate students
at Haverford College (``in an attempt to make the subject
accessible to everyone" \cite{C-F}). It is a pleasant method 
of coding representations of the fundamental
group of a knot complement into the group of symmetries 
of a regular triangle, but this interpretation is not needed 
in the definition and most of applications of 3-colorings.

\begin{definition}\label{1.15} (Fox 3-coloring of a link diagram).\\
Consider a coloring of a link diagram using colors r (red), 
y (yellow) or b (blue) 
in such a way that an arc of the diagram (from 
a tunnel to a tunnel) is colored by one color and at 
a crossing one uses one or all three colors. Such a coloring is
called a Fox 3-coloring. If whole diagram is colored by  just one 
color we say that we have a trivial coloring. Let $tri(D)$ denote 
the number of different Fox 3-colorings of $D$. 
\end{definition}
\begin{example}
\begin{enumerate}
\item[(i)] $tri(
\psfig{figure=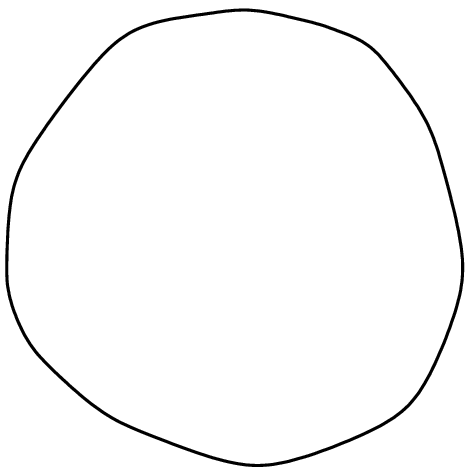,height=0.5cm}
) = 3$ as the trivial diagram has only trivial colorings.
\item[(ii)] $tri(
\psfig{figure=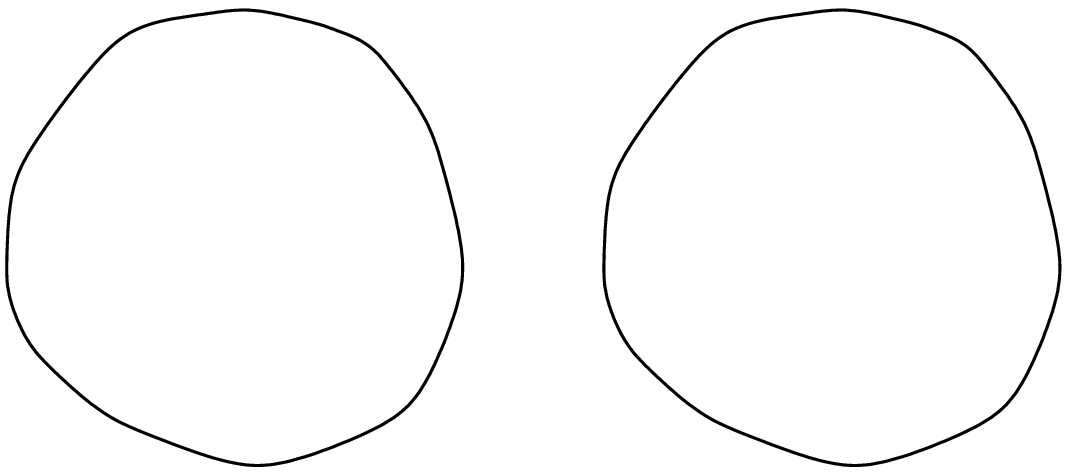,height=0.5cm}
)=9$, and more generally for a trivial link diagram of
$n$ components, $U_n$, one has $tri(U_n) = 3^n$.
\item[(iii)] For a standard diagram of a trefoil knot we have 
three trivial colorings
and 6 nontrivial colorings, one of them is presented in Figure 1.17 
(all other differ
from this one by permutations of colors. Thus tri(
\psfig{figure=aatrefoil.eps,height=.4cm})
$ = 3+6=9$
\end{enumerate}
\end{example}
\centerline{\psfig{figure=trefoil3-col.eps,height=3.1cm}}
\centerline{Fig. 1.17; \ Different colors are marked by lines of 
different thickness.}

Fox 3-colorings were defined for link diagrams, they are however invariants 
of links. One needs only show that $tri(D)$ 
is unchanged by Reidemeister moves.

The invariance under $R_1$ and $R_2$ is illustrated in Fig.1.18 and the
invariance under $R_3$ is illustrated in Fig.1.19.
\ \\
\ \\
\centerline{\psfig{figure=R1R23-col.eps}}
\begin{center}
Fig. 1.18
\end{center}

\centerline{\psfig{figure=R33-col.eps}}
\begin{center}
Fig. 1.19
\end{center}

The next property of Fox 3-colorings is the key in proving that
different trivial links are not 3-move equivalent.

\begin{lemma}[\cite{P-1}]
$tri(D)$ is unchanged by a 3-move.
\end{lemma}
The proof of the lemma is illustrated in Figure 1.20.

\ \\
\centerline{\psfig{figure=3-movecol.eps}}
\begin{center}
Fig. 1.20
\end{center}
The lemma also explains the fact that the trefoil has nontrivial
Fox 3-colorings: the trefoil knot is 3-move equivalent to
the trivial link of two components (Example 1.4(i)). 

Tomorrow I will place the theory of Fox coloring in more
general (sophisticated context), and apply it to the analysis 
of 3-moves (and $(2,2)$ and $(2,3)$ moves) of $n$-tangles.
Interpretation of tangle colorings as Lagrangians in symplectic
spaces is our main (and new) tool. In the second lecture tomorrow,
I will discuss another motivation for study 3-moves: 
to understand
skein modules based on their deformation.

\section{Talk 2: \ Lagrangian approximation of Fox $p$-colorings of tangles.}

We just have heard beautiful and elementary talk by
Lou Kauffman. I hope to follow his example by having my talk
elementary and deep at the same time. I will use several
results introduced by Lou, like rational tangles and their
classification, and will also
build on my yesterday's talk. The talk will culminate in the
introduction of the symplectic structure on the boundary 
of a tangle in such a way that tangles yields Lagrangians 
in the symplectic space. I could not dream of this
connection a year ago; however now, after 10 month perspective,
I see the symplectic structure as a natural development.

Let us start slowly from my personal perspective and motivation.
In the spring of 1986, I was analyzing behavior of Jones type invariants
of links when a link was modified by a $k$-move (or $t_k$, 
$\bar t_{2k}$ moves in the oriented case). My interest had its roots in the
fundamental Conway's paper \cite{Co}. In July of 1986, I gave a talk at
the "Braids" conference in Santa Cruz. I was told by Murasugi and Kawauchi,
after my talk, about the Nakanishi's 3-move conjecture.
It was suggested to me by R.Campbell (Kirby's student in 1986) 
to consider the effect of 3-moves on Fox colorings. 
Only several years later, when writting \cite{P-3} in 1993 
I realized that Fox colorings can be succesfully used to analyse 
moves on tangles, by considering not only the space of colorings 
but also the resulting coloring of boundary points. 
More of this later, now it is time to define Fox $k$-colorings.
\begin{definition}\label{2.1}
\begin{enumerate}
\item [(i)]
We say that a link (or tangle) diagram is k-colored if every
arc is colored by one of the numbers $0,1,...,k-1$ (forming a
group $Z_k$) in such a way that
at each crossing the sum of the colors of the undercrossings is equal
to twice the color of the overcrossing modulo $k$; see Fig.2.1. 
\item [(ii)]
The set of $k$-colorings form an abelian group, denoted by $Col_k(D)$.
The cardinality of the group will be denoted by $col_k(D)$.
For an $n$-tangle $T$ each Fox $k$-coloring of $T$ yields a
coloring of boundary points of $T$ and we have the homomorphism
$\psi :Col_k(T) \rightarrow Z_k^{2n}$
\end{enumerate}
\end{definition}.
\centerline{\psfig{figure=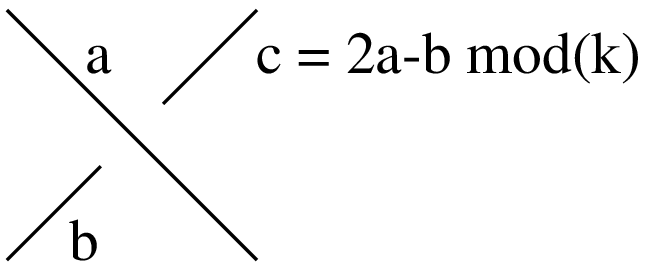,height=2.6cm}}
\centerline{Fig. 2.1}
It is a pleasant exercise to show that $Col_k(D)$ is unchanged
by Reidemeister moves and by $k$-moves (Fig.2.2). 
\\
\ \\
\centerline{\psfig{figure=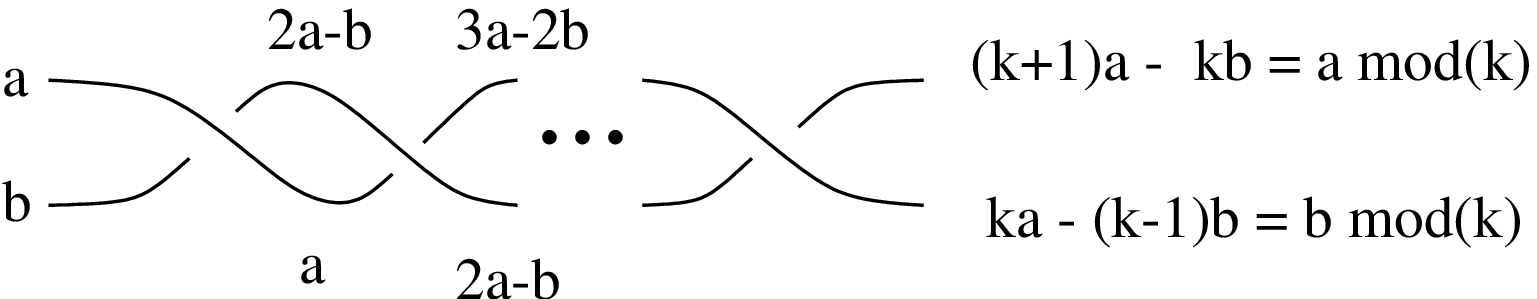,height=2.6cm}}
\centerline{Fig. 2.2}

Let us look more carefully at the observation that a $k$-move
preserves the space of Fox $k$-colorings and at the
unlinking conjectures described till now.
We discussed the 3-move conjecture and 
Nakanishi suggested also 4-move conjecture for knots (it does not hold
for links\footnote{Kawauchi suggested that for links one should conjecture
that two links are 4-move equivalent iff they are link homotopic.}).
As I mentioned yesterday not every link can be simplified using
$5$-moves, but a $5$-move is a combination of $(2,2)$ moves and these
moves possibly suffices to reduce every link. Similarly not
every link can be reduced using
$7$-moves, but a $7$-move is a combination of 
$(2,3)$-moves\footnote{To be precise, a $7$-move is a combination
of a $(-3,-2)$ and $(2,3)$ moves; compare Fig.1.10.}
which possibly suffice for reduction. We stopped at this
point yesterday, but what can one use instead of a general $k$-move?
Let us consider the case of a $p$-move where $p$ is a prime number.
I suggest (and state publicly for the first time)
that one should consider {\it rational moves}, that is,
a rational $\frac{p}{q}$-tangle of Conway is substituted in
place of the identity tangle\footnote{The move was first
considered by Y.Uchida \cite{Uch}.}. 
The important observation for us is that
$Col_{p}(D)$ is preserved by $\frac{p}{q}$-moves. Fig.2.3
illustrates the fact that $Col_{13}(D)$ is unchanged by
a $\frac{13}{5}$-move.
\\
\ \\
\centerline{\psfig{figure=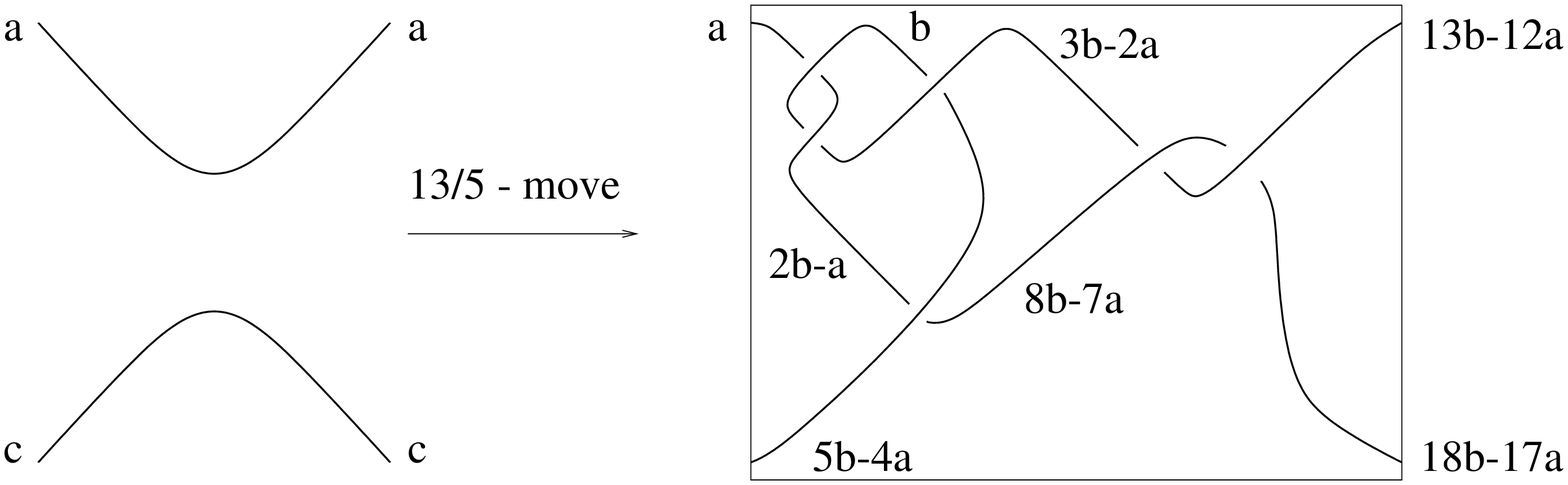,height=3.6cm}}
\centerline{Fig. 2.3}

We should note also that $(m,q)$-moves are equivalent
to $\frac{mq+1}{q}$-moves (Fig.2.4) 
so the space of Fox $(mq+1)$-colorings
is preserved.
\\
\ \\
\centerline{\psfig{figure=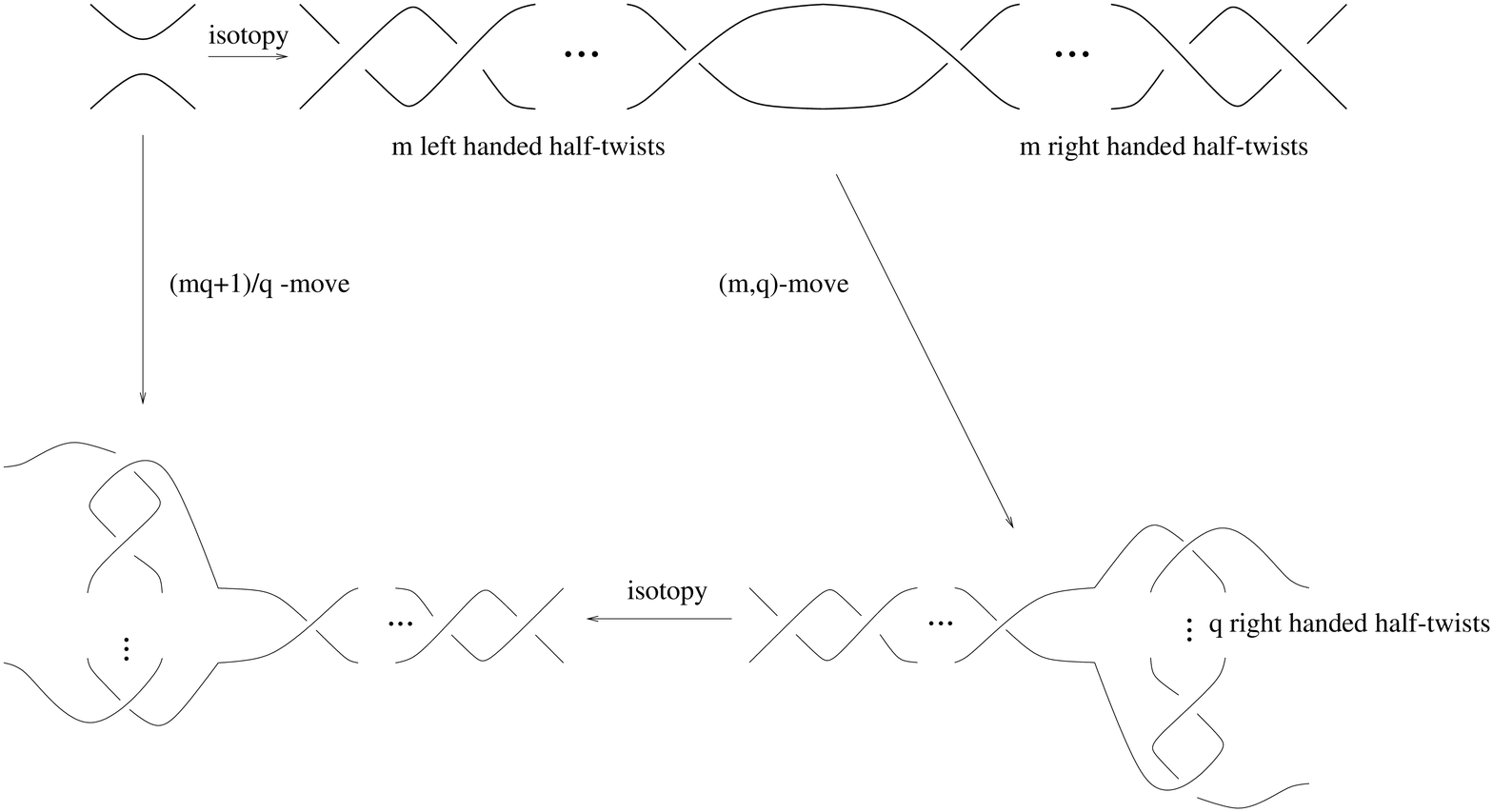,height=4.9cm}}
\centerline{Fig. 2.4}

We just have heard about the Conway's  classification
of rational tangles at the Lou's talk, so  I will
briefly sketch definitions and notation.
The 2-tangle of Figure 2.5 is called a rational tangle
with Conway notation $T(a_1,a_2,...,a_n)$. It is a rational
$\frac{p}{q}$-tangle if $\frac{p}{q} = 
a_n + \frac{1}{a_{n-1}+...+\frac{1}{a_1}}$.\footnote{$\frac{p}{q}$
is called the slope of the tangle and can be easily
identified with the slope of the meridian disk of the solid torus
being the branched double cover of the rational tangle.} Conway
proved that two rational tangles are ambient isotopic
(with boundary fixed) if and only if their slopes
are equal (compare \cite{Kaw}).
\\
\ \\
\centerline{\psfig{figure=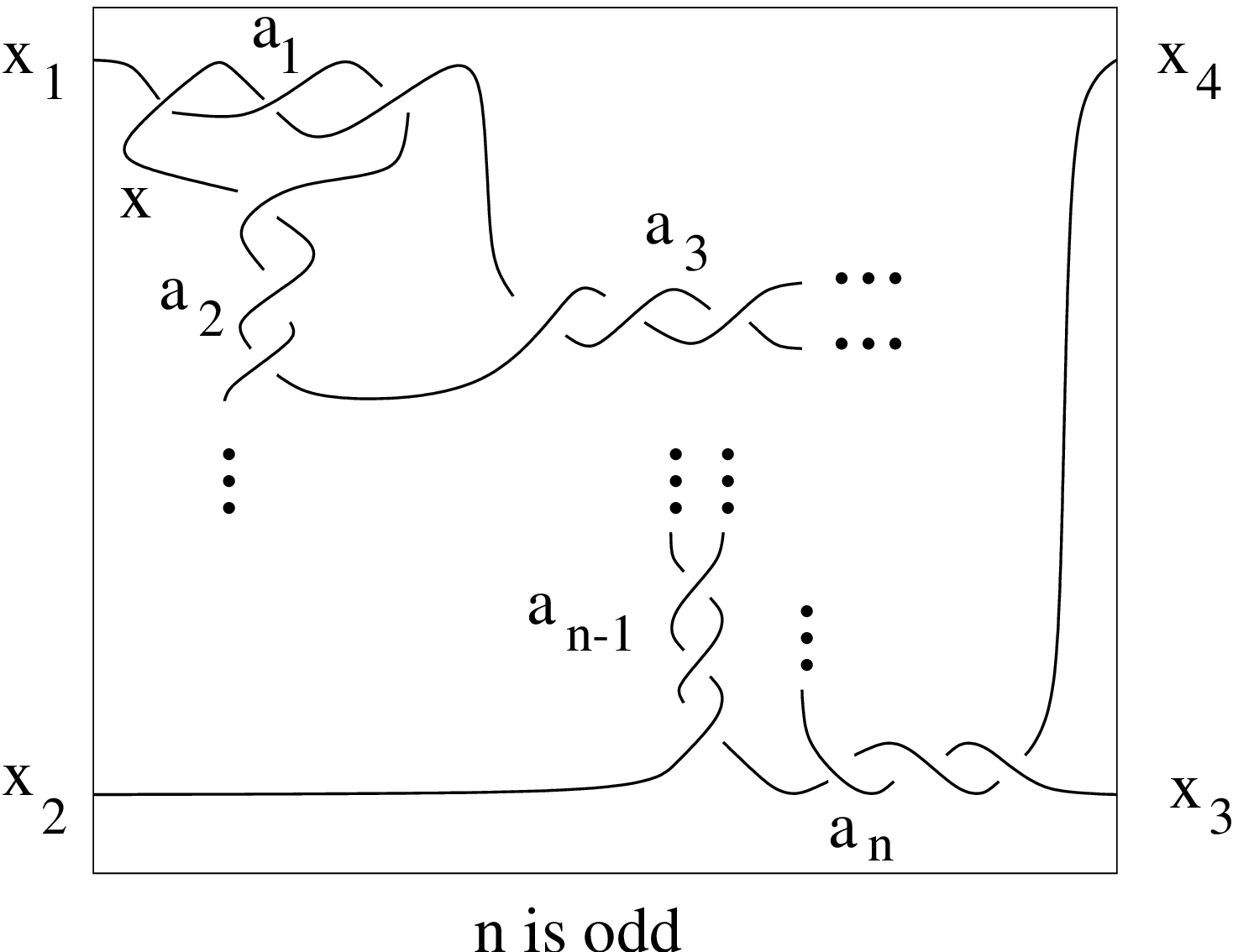,height=4.6cm}\ \ \ 
\psfig{figure=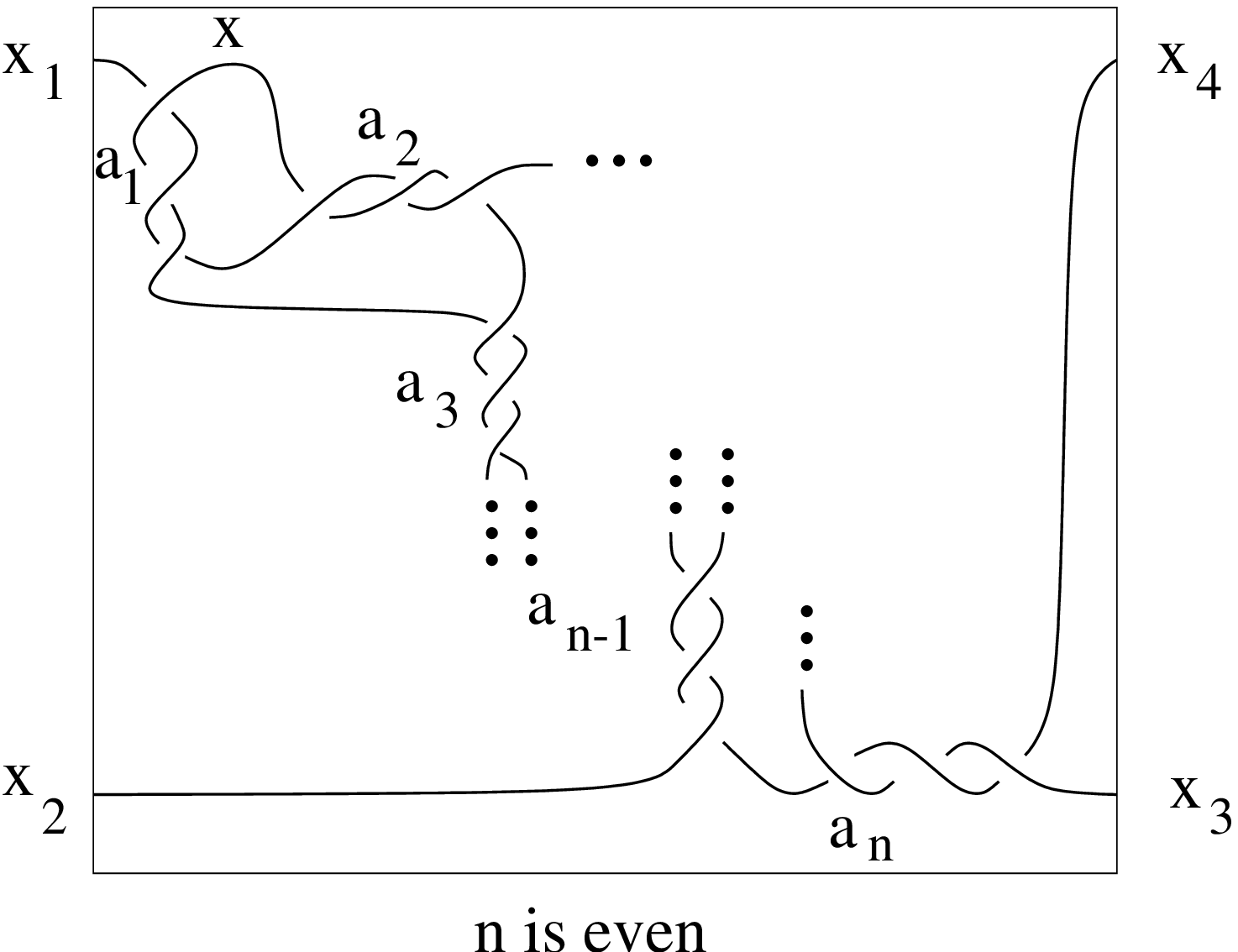,height=4.9cm}}
\centerline{Fig. 2.5}

For a Fox coloring of a rational $\frac{p}{q}$-tangle with
boundary colors $x_1,x_2,x_3,x_4$ (Fig.2.5), one has 
$x_4-x_1 = p(x-x_1)$, $x_2-x_1=q(x-x_1)$ and $x_3 = 
x_2 + x_4 -x_1$. If a coloring is nontrivial ($x_1\neq x$) 
then $\frac{x_4-x_1}{x_2-x_1} = \frac{p}{q}$ as explained by
Lou.

\begin{conjecture}\label{2.2}\ \\
Let $p$ be a fixed prime number, then
\begin{enumerate} 
\item[(i)] Every link can
be reduced to a trivial link
by rational $\frac{p}{q}$-moves ($q$ any integer).
\item[(ii)] There is a function
$f(n,p)$ such that any $n$-tangle can be reduced to one of
``basic" $f(n,p)$ $n$-tangles (allowing additional trivial components)
by rational $\frac{p}{q}$-moves.
\end{enumerate}
\end{conjecture}
First we can observe that it suffices to consider $\frac{p}{q}$-moves with
$|q| \leq \frac{p}{2}$, as other $\frac{p}{q}$-moves follow. Then we know
that for $p$ odd the $\frac{p}{1}$-move is a combination of $\frac{p}{2}$ and
$\frac{p}{-2}$-moves (Compare Fig.1.8). Thus, in fact, 3-move, 
$(2,2)$-move and
$(2,3)$-move conjectures are special cases of Conjecture 2.2(i).
If we analyze the case $p=11$ we see that $\frac{11}{2} = 5+ \frac{1}{2}$,
$\frac{11}{3} = 4 - \frac{1}{3}$, $\frac{11}{4} = 3 - \frac{1}{4}$,
$\frac{11}{5} = 2+ \frac{1}{5}$. Thus:
\begin{conjecture}\label{2.3}\ \\
Every link can be reduced to a trivial link (with the same 
space of 11-colorings)
by $(2,5)$ and $(4,-3)$ moves, they inverses and mirror images.
\end{conjecture}

What about the number $f(n,p)$?
We know that because $\frac{p}{q}$-moves preserve $p$-colorings, 
therefore $f(n,p)$ is bounded
from below by the number of subspaces of $p$-colorings 
of the $2n$ boundary points induced by Fox $p$-colorings 
of $n$-tangles (that is by the number
of subspaces $\psi (Col_p(T))$ in $Z_p^{2n}$).
I noted in \cite{P-3} that for 2-tangles this number is equal to $p+1$
(even in this special case my argument was complicated).
For $p=3,n=4$
the number of subspaced followed from the work of my student
Tatsuya Tsukamoto and was equal to $40$ \cite{P-Ts}.
The combined effort of Mietek D{\c a}bkowski and Tsukamoto
gave the number $1120$ for subspaces $\psi (Col_3(T))$ and 4-tangles.
This was my knowledge at the early spring of 2000.
On May 2 and 3 I heard talks on Tits buildings
(at the Banach Center in Warsaw) by J.Dymara and T.Januszkiewicz. 
I realized that the topic may have some 
connection to my work.
I asked Tadek Januszkiewicz whether he sees relations 
and I gave him numbers $4$, $40$, $1120$ for $p=3$.
He immediately answered that most likely I counted the number
of Lagrangians in $Z_3^{2n-2}$ symplectic space, and that
the number of Lagrangians in $Z_p^{2n-2}$ is known to be
equal to $\prod_{i=1}^{n-1} (p^i +1)$.
Soon I constructed the appriopiate symplectic form (as did
Janek Dymara). I will spend most of this talk on the
construction. I will end with discussion of classes
of tangles for which it has been proven that 
$f(n,p) = \prod_{i=1}^{n-1} (p^i +1)$.

 Consider $2n$ points on a circle (or a square) and a
field ${\bf Z}_p$ of $p$-colorings of a point. Thus colorings of
$2n$ points form ${\bf Z}_{p}^{2n}$ linear space. 
Let $e_1,\ldots , e_{2n}$ be its basis,
\ \\
\centerline{\psfig{figure=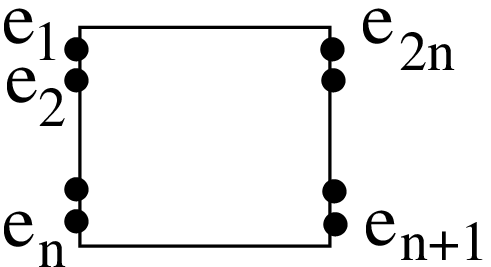}}
\begin{center} Fig. 2.6 \end{center}
\ \\
 $e_i=(0,\ldots , 1,\ldots ,0)$, where 1
occurs in the $i$-th position. Let ${\bf Z}_{p}^{2n-1}\subset{\bf
Z}_p^{2n}$ be the subspace of vectors $\sum a_i e_i$ satisfying $\sum (-1)^i
a_i=0$ (alternating condition). 
Consider the basis $f_1,\ldots , f_{2n-1}$ of ${\bf
Z}_{p}^{2n-1}$ where $f_k=e_k+e_{k+1}$. Consider a skew-symmetric
form $\phi$ on ${\bf Z}_{p}^{2n-1}$ of nullity 1 given by the
matrix

$$\phi = \left( \begin{array}{cccc} 0 & 1 & \ldots &\ldots \\
 -1 & 0 & 1 &\ldots \\ \ldots & \ldots & \ldots & \ldots \\ \ldots & \ldots
 & -1 & 0
\end{array} \right)$$

that is
 \vspace{1mm}
   \renewcommand{\arraystretch}{3}
   $$\phi (f_i, f_j) =\left\{
   \begin{array}{lr}
   0 &
   {\rm if}\ |j-i|\neq 1 \\
  1 &  {\rm if\ }\ j=i+1\\ -1 &
   {\rm if}\ j=i-1. \\
   \end{array}
  \right .$$
   \par
   \vspace{2mm}

Notice that the vector $e_1+ e_2+ \ldots + e_{2n}$
($=f_1+f_3+\ldots + f_{2n-1}=f_2+f_4+\ldots + f_{2n})$
is $\phi$-orthogonal to any other vector.
If we consider ${\bf Z}_{p}^{2n-2}={\bf Z}_{p}^{2n-1}/{\bf
Z}_{p}$, where the subspace ${\bf Z}_{p}$ is generated by
$e_1+\ldots + e_{2n}$, that is, ${\bf Z}_{p}$ consists
of monochromatic (i.e. trivial) colorings,
then $\phi$ descends to the symplectic form $\hat\phi$ on ${\bf
Z}_{p}^{2n-2}$. Now we can analyze isotropic subspaces of
$({\bf Z}_{p}^{2n-2},\hat\phi)$, that is subspaces on
which  $\hat\phi$ is $0$
($W\subset {\bf Z}_{p}^{2n-2}, \phi (w_1,w_2)=0$ for $w_1,w_2\in W$).
The maximal isotropic ($(n-1)$-dimensional) 
subspaces of ${\bf Z}_{p}^{2n-2}$ are
called Lagrangian subspaces (or maximal totally degenerated subspaces)
and there are $\prod_{i=1}^{n-1}(p^i+1)$ of them.

We have $\psi :Col_p T\rightarrow{\bf Z}_{p}^{2n}$. Our local
condition on Fox colorings (Fig.2.1) guarantees that for any tangle $T$,
$\psi (Col_p T)\subset {\bf Z}_{p}^{2n-1}$. Also
${\bf Z}_{p}$, the space of trivial colorings, always lays in $Col_p T$.
Thus $\psi$ descents to $\hat\psi :Col_p T/{\bf Z}_{p} \rightarrow
{\bf Z}_{p}^{2n-2}={\bf Z}_{p}^{2n-1}/{\bf Z}_{p}$ Now we have
the fundamental question:
Which subspaces of
${\bf Z}_{p}^{2n-2}$ are yielded by $n$-tangles?
We answer this question below.
\begin{theorem}
$\hat\psi (Col_p T/{\bf Z}_{p})$ is a Lagrangian
subspace of ${\bf Z}_{p}^{2n-2}$ with the symplectic form $\hat\psi$.
\end{theorem}
The natural question would be whether every Lagrangian subspace
can be realized by a tangle. The answer is negative for
$p=2$ and positive for $p>2$ \cite{D-J-P}.
As a corollary we obtain a fact which before was considered to
be difficult even for 2-tangles.
\begin{corollary}
For any $p$-coloring of a tangle boundary satisfying 
the alternating property (i.e. an element
of ${\bf Z}_{p}^{2n-1}$) there is an $n$-tangle and
its $p$-coloring yielding the given coloring on the boundary.
In other worlds: ${\bf Z}_{p}^{2n-1} = \bigcup_T \psi_T(Col_p(T))$.
Furthermore the space $\psi_T(Col_p(T))$ is $n$-dimensional.
\end{corollary}
We can say that we understand the conjectured value of 
the function $f(n,p)$ but when can we prove 
Conjecture 2.2 with $f(n,p)= \prod_{i=1}^{n-1}(p^i+1)$?
In fact we know that for $p=2$ not every Lagrangian 
is realized and actually
$f(n,2)= \prod_{i=1}^{n-1}(2i+1)$. 
For rational 2-tangles Conjecture 2.2 follows almost from the definition and
the generalization to 2-algebraic tangles (algebraic tangles in the Conway
sense) is not difficult. In order to be able to use induction for 
tangles with $n>2$ we generalize the concept of the algebraic tangle:
\begin{definition}\label{2.4}\ \\
\begin{enumerate}
\item[(i)] n-algebraic tangles is the smallest family 
of n-tangles which satisfies:\\
(0) Any n-tangle with 0 or 1 crossing is n-algebraic.\\
(1) If $A$ and $B$ are n-algebraic tangles then $r^i(A)*r^j(B)$ 
is n-algebraic; $r$ denotes here the rotation of a tangle by 
$\frac{2\pi}{2n}$ angle, and * denotes (horizontal) composition of tangles.

\item[(ii)] If in the condition (1), $B$ is restricted 
to tangles with no more than $k$ crossings, we obtain 
the family of $(n,k)$-algebraic tangles.

\item[(iii)] 
If an $m$-tangle, $T$, is obtained from an $(n,k)$-algebraic 
tangle (resp. $n$-algebraic tangle)
by partially closing its endpoints ($2n-2m$ of them) 
without introducing any new crossings then $T$ is called
an $(n,k)$-algebraic (resp. $k$-algebraic) $m$-tangle. For $m=0$
we obtain an $(n,k)$-algebraic (resp. $k$-algebraic) link.
\end{enumerate}
\end{definition}
Conjecture 2.2, for $p=3$, has been proven for 3-algebraic 
tangles \cite{P-Ts}  ($f(3,3)=40$) and $(4,5)$-algebraic 
tangles \cite{Tsu} ($f(4,3)=1120$). In particular the 
Montesinos-Nakanishi conjecture holds for 3-algebraic and
$(4,5)$-algebraic links. 40 ``basic" 3-tangles are
shown in Fig. 2.7.\\
\ \\
\centerline{\psfig{figure=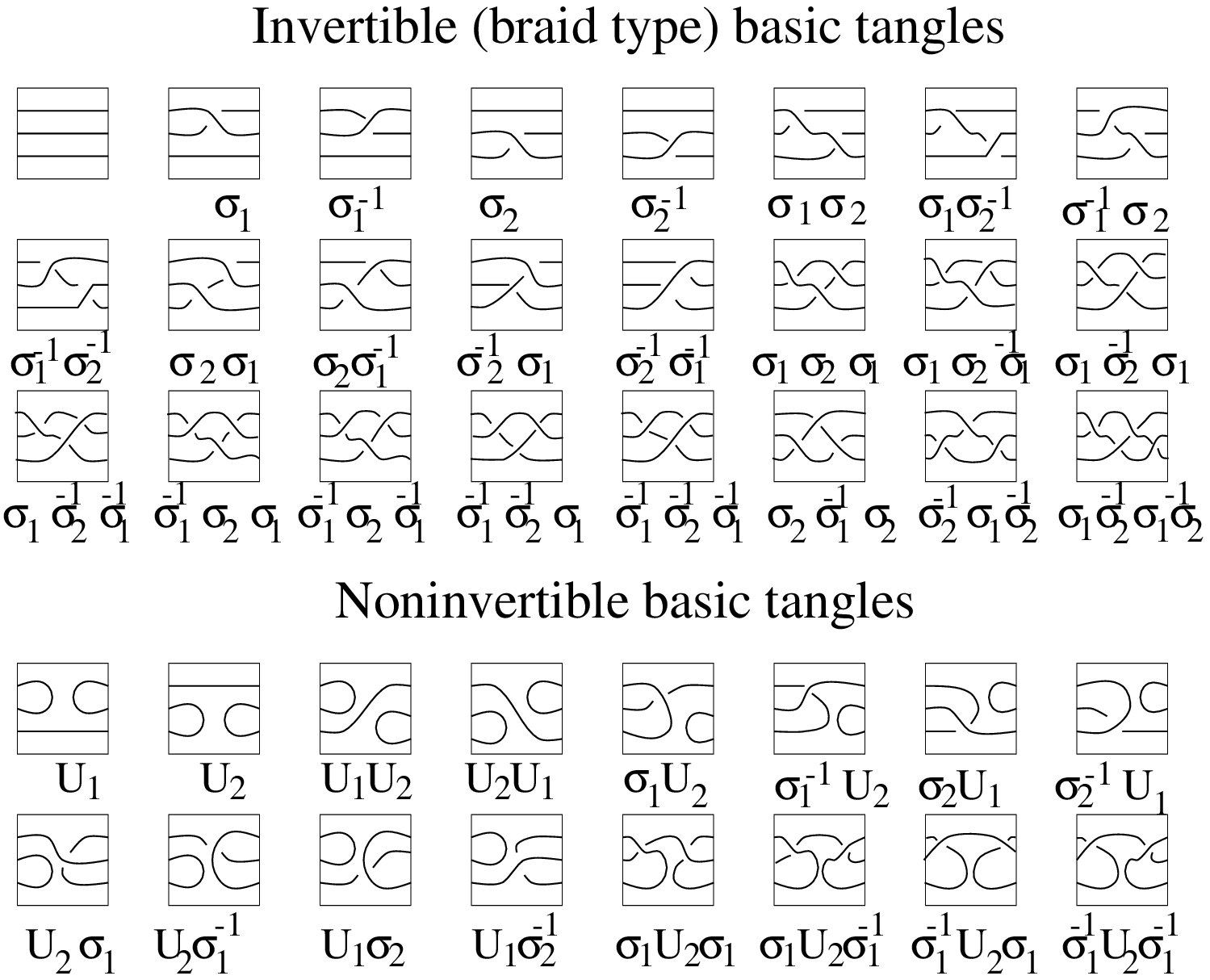,height=8.8cm}}
\centerline{Fig. 2.7}

The simplest $4$-tangles which cannot be distinguished by
3-coloring for which 3-move equivalence is not yet established
are illustrated in Fig.2.8. With respect to $(2,2)$ moves,
equivalence of 2-tangles of Fig.2.9 is still an open problem.

\centerline{\psfig{figure=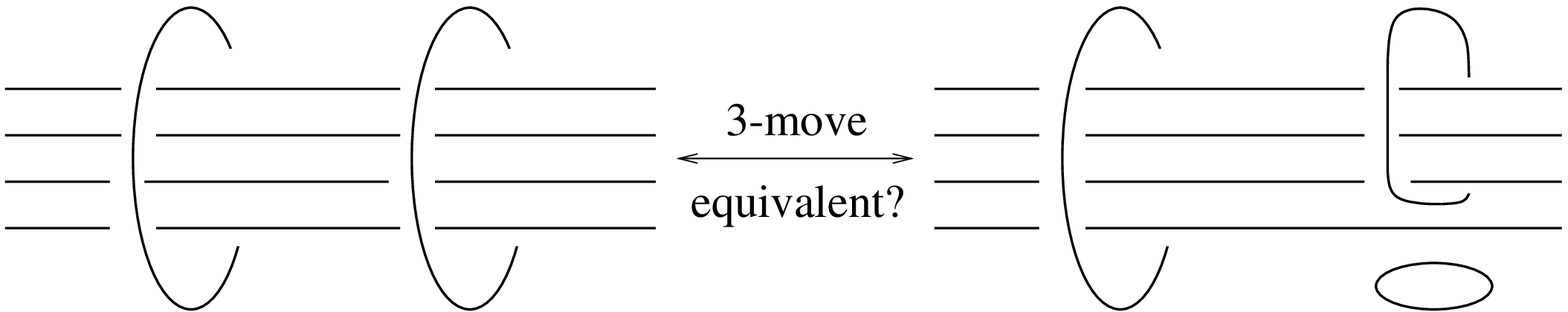,height=3.2cm}}
\centerline{Fig. 2.8}
\centerline{\psfig{figure=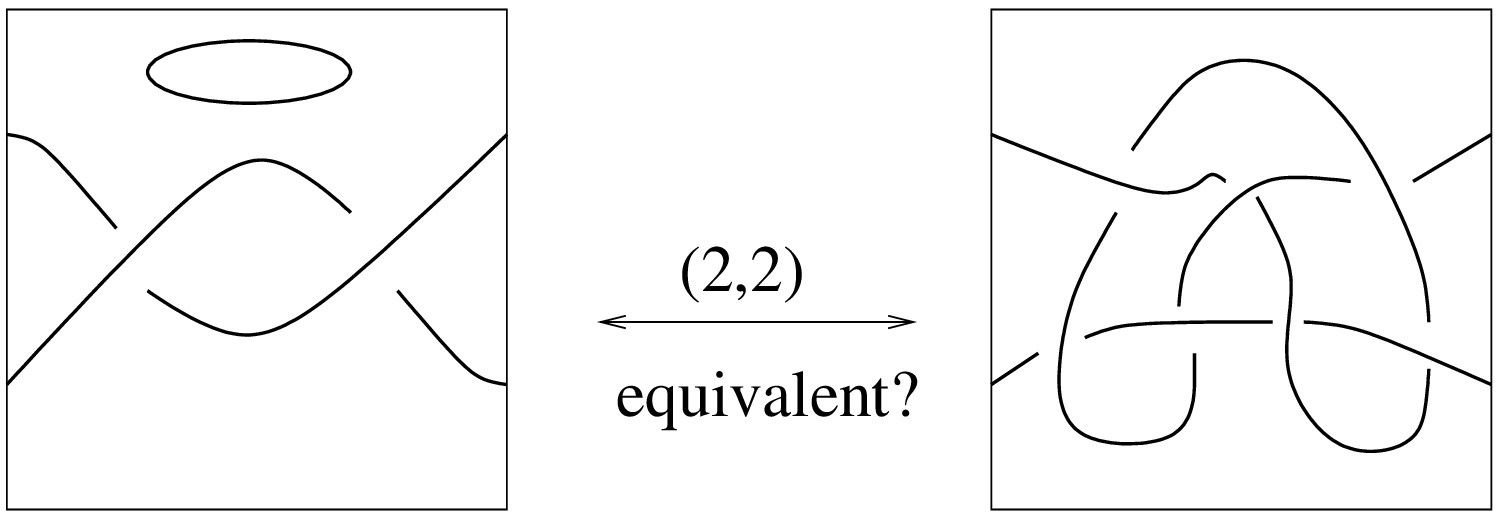,height=3.8cm}}
\ \\
\ \\
\centerline{Fig. 2.9}

\section{Talk 3: \ Historical Introduction to Skein Modules.}

I will discuss, in my last talk of the conference,
{\it skein modules}, or as I prefer to say more generally, 
{\it algebraic topology based on knots}. It is my 
``brain child" even if the idea was also conceived by
other people (most notably Vladimir Turaev), and was
envisioned by John H.Conway (as ``linear skein") a
decade earlier.
Skein modules have their origin in the observation by Alexander
\cite{Al}, that his polynomials ({\it Alexander polynomials}) of three
links, $L_+,L_-$ and $L_0$ in $R^3$ are linearly related (Fig.3.1).

For me it started in Norman Oklahoma in April of 1987,
when I was enlightened to see that the multivariable
version of the Jones-Conway (Homflypt) polynomial
analyzed by Hoste and Kidwell is really a module
of links in a solid torus (or more generally in
the connected sum of solid tori).

I would like to discuss today, in more detail,  
skein modules related to the (deformations) 
of 3-moves and the Montesinos-Nakanishi
conjecture but first I will give the general definition
and I will make a short tour of the world of skein modules.

Skein Module is an algebraic object associated to a manifold, 
usually constructed as
a formal linear combination of embedded (or immersed) submanifolds,
modulo locally defined relations.
In a more restricted setting a {\it skein module}
is a module associated to a 3-dimensional manifold, by considering
linear combinations of links in the manifold, modulo properly chosen
(skein) relations. It is the main object of the {\it algebraic topology
based on knots}. In the choice of relations one takes
into account several factors:
\begin{enumerate}
\item [(i)]  Is the module we obtain accessible (computable)?
\item [(ii)]
How precise are our modules in distinguishing 3-manifolds and links
in them?
\item [(iii)] Does the module reflect topology/geometry of a 3-manifold
(e.g. surfaces in a manifold, geometric decomposition of a manifold)?.
\item [(iv)]
Does the module admit some additional structure (e.g.~filtration,
gradation, multiplication, Hopf algebra structure)? Is it leading
to a Topological Quantum Field Theory (TQFT) by taking a finite quotient?
\end{enumerate}
\ \\
One of the simplest skein modules is a $q$-deformation of the first
homology group of a 3-manifold $M$, denoted by ${\cal S}_2(M;q)$.
It is based on the skein relation (between unoriented framed links in $M$):
$L_+=qL_0$; it also satisfies the framing relation 
$L^{(1)}-qL$, where $L^{(1)}$
denote a link obtained from  $L$ by twisting the framing of $L$ once
in the positive direction. 
Already this simply defined skein module ``sees"
nonseparating surfaces in $M$. These surfaces are responsible for
torsion part of our skein module \cite{P-10}.

There is more general pattern: most of analyzed skein modules reflect
various surfaces in a manifold.

The best studied skein modules use skein relations
which worked successfully in the classical
knot theory (when defining polynomial invariants of links in $R^3$).
\begin{enumerate}
\item[(1)] The Kauffman bracket skein module, KBSM.\\
The skein module based on the {\it Kauffman bracket skein relation},
$L_+ = AL_- + A^{-1}L_{\infty}$, and denoted by $S_{2,\infty}(M)$,
is best understood among the Jones type skein modules.
It can be interpreted as a quantization of the co-ordinate ring
of the character variety of $SL(2,C)$ representations of the
fundamental group of the manifold $M$, \cite{Bu-2,B-F-K,P-S}.
For $M= F\times [0,1]$,
KBSM is an algebra (usually noncommutative). It is finitely generated
algebra for a compact $F$ \cite{Bu-1}, and has no zero divisors \cite{P-S}.
The center of the algebra is generated by boundary 
components of $F$ \cite{B-P,P-S}.
Incompressible tori and 2-spheres in $M$ yield torsion in KBSM;
it is a question of fundamental importance whether other surfaces
can yield torsion as well.

\item[(2)] Skein modules based on the Jones-Conway (Homflypt) relation.\\
$v^{-1}L_+ - vL_- = z L_0$, where $L_+,L_-,L_0$ are oriented links 
(Fig.3.1).
These skein modules are denoted by $S_3(M)$ and generalize skein modules
based on Conway relation which were hinted by Conway. For $M= F\times [0,1]$,
$S_3(M)$ is a {\it Hopf algebra} (usually neither
commutative nor co-commutative),
\cite {Tu-2,P-6}. $S_3(F\times [0,1])$ is a free module and can be
interpreted as a quantization \cite{H-K,Tu-1,P-5,Tu-2}.
$S_3(M)$ is related to the algebraic set of
$SL(n,C)$ representations of the
fundamental group of the manifold $M$, \cite{Si}.\\
\ \\
\centerline{\psfig{figure=L+L-L0.eps}}
\centerline{Fig. 3.1}
\item[(3)] Skein modules based on the {\it Kauffman polynomial} relation\\
$L_{+1}  + L_{-1} = x (L_0 + L_{\infty})$ (see Fig.3.3) and the framing 
relation $L^{(1)}-aL$.
It is denoted by $S_{3,\infty}$ and is known to be free for $M= F\times [0,1]$.

\item[(4)] Homotopy skein modules. In these skein modules, $L_+ = L_-$ for
self-crossings. The best studied example is the q-homotopy skein module
with the skein relation $q^{-1}L_+ -qL_- =zL_0$ for mixed crossings.
For $M= F\times [0,1]$ it is a quantization, \cite{H-P-1,Tu-2,P-11}, and
as noted by Kaiser they can be almost completely understood using
singular tori technique of Lin.

\item[(5)]
Skein modules based on Vassiliev-Gusarov filtration.\\
We extend the family of knots, $\cal K$, by singular knots,
and resolve s singular crossing by $K_{cr}= K_+ -K_-$. These allows us
to define the Vassiliev-Gusarov filtration:
$... \subset C_3  \subset C_2  \subset C_1  \subset C_0= R{\cal K}$,
where $C_k$ is generated by knots with $k$ singular points.
The k'th Vassiliev-Gusarov skein module is defined to be a quotient:
$W_k(M) = R{\cal K}/C_{k+1}$. The completion of the space of knots
with respect to the Vassiliev-Gusarov filtration, $\hat{R{\cal K}}$,
 is a {\it Hopf algebra} (for $M=S^3$). Functions dual to
Vassiliev-Gusarov skein modules are called {\it finite type} or
{\it Vassiliev invariants} of knots; \cite{P-7}.

\item[(6)]
Skein modules based on relations deforming n-moves.\\
${\cal S}_n(M)= R{\cal L}/(b_0L_0 + b_1L_1 + b_2L_2 +...+b_{n-1}L_{n-1})$.
In the unoriented case, we can add to the relation the 
term $b_{\infty}L_{\infty}$,
to get ${\cal S}_{n,\infty}(M)$, and also,
possibly, a framing relation.  The case $n=4$, on which I am working
with my students,  will be described, n greater detail, 
in a moment.
\end{enumerate}
\
\\
Examples (1)-(5) gave a short description of skein 
modules studied extensively untill now. I will now spent 
more time on two other examples which only recently
has been considered more detaily. The first example is based on
a deformation of the 3-move and the second on the deformation of
the $(2,2)$-move. The first examples has been studied by
my students Tsukamoto and Mike Veve. I denote it by 
${\cal S}_{4,\infty}$ as it involves (in the skein relation),
4 horizontal positions and the vertical ($\infty$) smoothing.
\begin{definition}
Let $M$ be an oriented 3-manifold, ${\cal L}_{fr}$ the set of unoriented
framed links in $M$ (including the empty knot, $\emptyset$) and
$R$ any commutative ring with identity.  Then we define the
$(4,\infty)$ skein module as:
${\cal S}_{4,\infty}(M;R) = R{\cal L}_{fr}/I_{(4,\infty)}$,
where $I_{(4,\infty})$ is the submodule of $R{\cal L}_{fr}$
generated by the skein relation:\\
$b_0L_0  + b_1L_1 + b_2L_2 + b_3L_3 + b_{\infty}L_{\infty} = 0$
and
the framing relation:\\
$L^{(1)} = a L$ where $a,b_0,b_3$ are invertible elements in $R$ and
$b_1,b_2,b_{\infty}$ are any fixed elements of $R$ (see Fig.3.2).
\end{definition}
\centerline{\psfig{figure=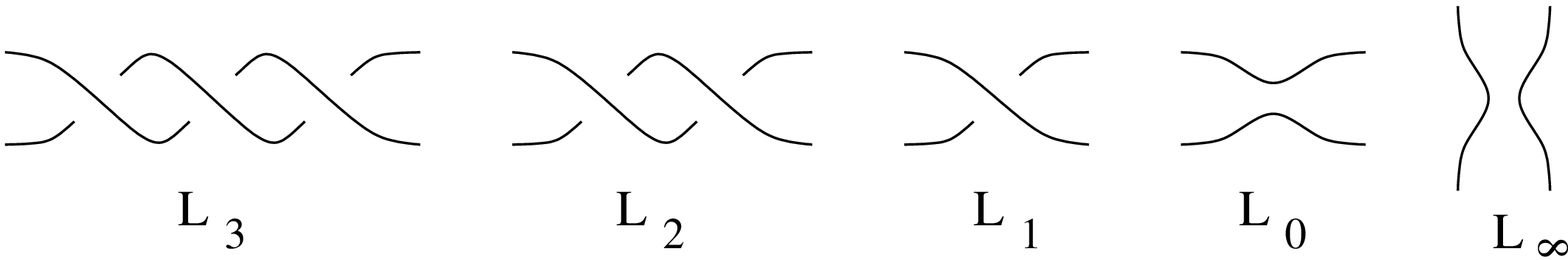,height=2.2cm}}
\centerline{Fig. 3.2}

The generalization of the Montesinos-Nakanishi cojnecture
says that ${\cal S}_{4,\infty}(S^3,R)$ is generated by
trivial links and that for $n$-tangles our skein
module is generated by $f(n,3)$ basic tangles (with
possible trivial components). This would give a
generating set for our skein module of $S^3$ or 
$D^3$ with
$2n$ boundary points (an $n$-tangle).
In \cite{P-Ts} we analyze extensively the possibilities
that trivial links, $T_n$, are linearly independent.
This may happen if $b_{\infty} = 0$ and $b_0b_1=b_2b_3$.
These leads to the following conjecture:
\begin{conjecture}\label{3.2}
\begin{enumerate}
\item [(1)] There is a polynomial invariant of unoriented links,
$P_1(L) \in Z[x,t]$ which satisfies:
\begin{enumerate}
\item [(i)] Initial conditions: $P_1(T_n) = t^n$, where $T_n$ is a trivial
link of $n$ components.
\item [(ii)] Skein relation $P_1(L_0) + xP_1(L_1) - xP_1(L_2) - P_1(L_3)=0$
where $L_0,L_1,L_2,L_3$ is a standard, unoriented skein quadruple
($L_{i+1}$ is obtained from $L_{i}$ by a right-handed half twist on
two arcs involved in $L_{i}$; compare Fig.3.2).
\end{enumerate}
\item [(2)] There is a polynomial invariant of unoriented framed links,
$P_2(L) \in Z[A^{\pm 1},t]$ which satisfies:
\begin{enumerate}
\item [(i)] Initial conditions: $P_2(T_n) = t^n$,
\item [(ii)] Framing relation: $P_2(L^{(1)})=-A^3P_2(L)$ where $L^{(1)}$ is
obtained from a framed link $L$ by a positive half twist on its framing.
\item [(iii)] Skein relation: $P_2(L_0) + A(A^2 + A^{-2})P_2(L_1) +
(A^2 + A^{-2})P_2(L_2) + AP_2(L_3)=0$.
\end{enumerate}
\end{enumerate}
\end{conjecture}
The above conjectures assume that $b_{\infty}=0$ in our
skein relation. Let us consider, for a moment, the possibility that
$b_{\infty}$ is invertible in $R$. Using the ``denominator"
of our skein relation (Fig.3.3) we get the relation which
allows us to compute the effect of adding a trivial component
to a link $L$ (we write $t^n$ for the trivial link link $T_n$):
$$(*)\ \ \ (a^{-3}b_3 + a^{-2}b_2 + a^{-1}b_1 + b_0 + b_{\infty}t)L=0$$
When considering the ``numerator" of the relation and its 
mirror image (Fig.3.3)
we obtain formulas for Hopf link summands, and because the unoriented
Hopf link is amphicheiral thus we can eliminate it from
our equations to get the formula (**):
$$b_3(L\#H) + (ab_2 + b_1t +a^{-1}b_0 + ab_{\infty})L =0.$$
$$b_0(L\#H) + (a^{-1}b_1 +b_2t + ab_3 + a^2b_{\infty})L =0.$$
$$(**)\ \ \ ((b_0b_1 - b_2b_3)t + (a^{-1}b_0^2 -a b_3^2) +
(ab_0b_2 - a^{-1}b_1b_3) + b_{\infty}(ab_0 -a^2b_3))L = 0.$$
\ \\
\centerline{\psfig{figure=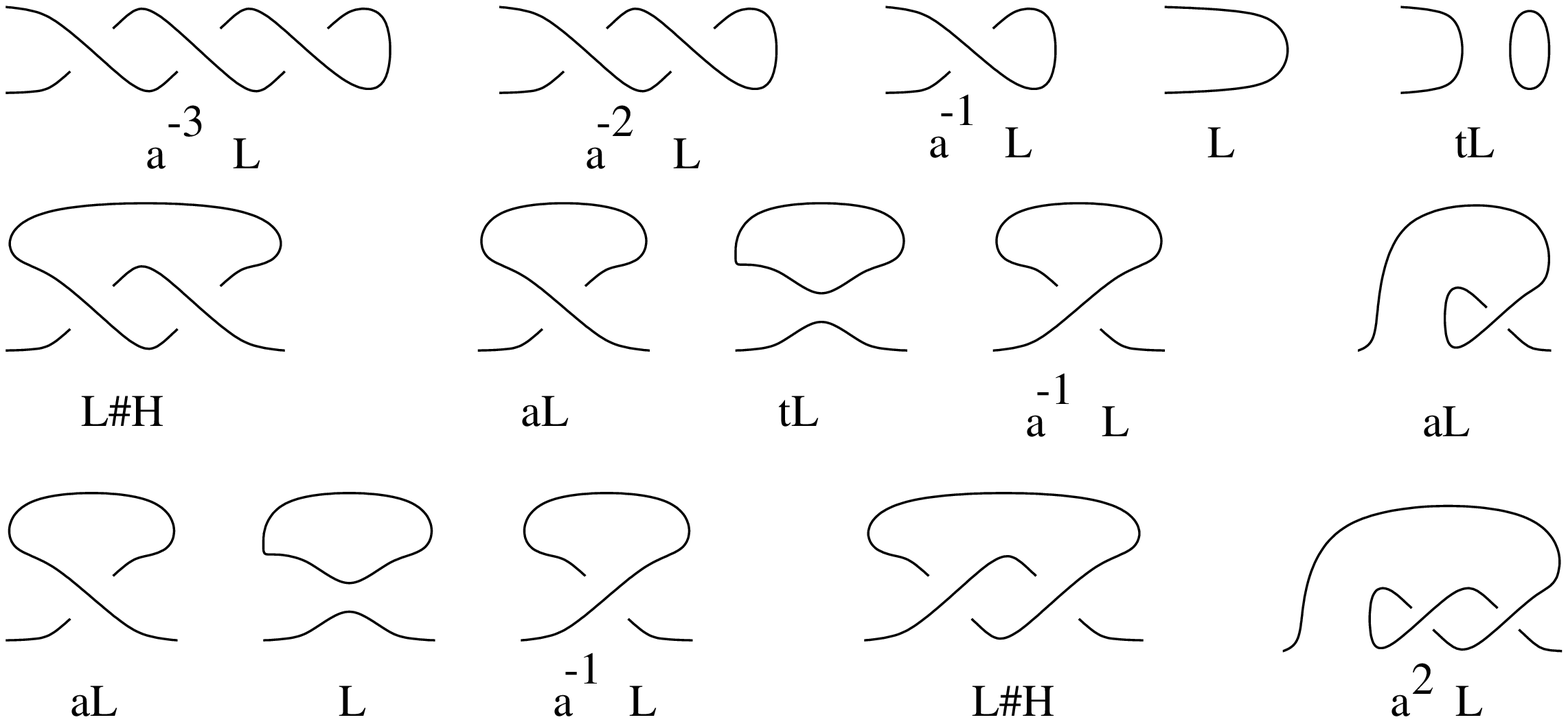,height=5.6cm}}
\centerline{Fig. 3.3}
It is possible that (*) and (**) are the only relations in the module.
The pleasant substitution which realizes the relation is:
$b_0=b_3=a=1$, $b_1=b_2=x$, $b_{\infty}=y$. This may lead to the
polynomial invariant of unoriented links in $S^3$ with value in $Z[x,y]$ 
and the skein relation $L_3 +xL_2 + xL_1 + L_0 + yL_{\infty} = 0$.

What about the relations to Fox colorings?
One such relation, already mentioned,
is the use of 3-coloring to estimate the number of
basic n-tangles (by $\prod_{i=1}^{n-1} (3^i +1)$) for
the skein module ${\cal S}_{4,\infty}$. I am also convinced that
${\cal S}_{4,\infty}(S^3;R)$ contains full information on the space of
Fox 7-colorings. It would be generalization of the fact that
the Kauffman bracket polynomial contains information on 3-colorings
and the Kauffman polynomial contains information on 
5-colorings. In fact, Fran{\c c}ois Jaeger told me that he 
knew how to get the space of $p$-colorings from
a short  skein relation (of type of $(\frac{p+1}{2},\infty)$).
Unfortunately
Fran{\c c}ois died prematurely in 1997 and I do not know 
how to prove his statement.

Finally let me describe shortly the skein module
related to the $(2,2)$-move conjecture.
Because a $(2,2)$-move is equivalent to the rational
$\frac{5}{2}$-move, I will denote the skein module
by ${\cal S}_{\frac{5}{2}}(M;R)$.
\begin{definition}
Let $M$ be an oriented 3-manifold, ${\cal L}_{fr}$ the set of unoriented
framed links in $M$ (including the empty knot, $\emptyset$) and
$R$ any commutative ring with identity.  Then we define the 
$\frac{5}{2}$-skein module as:
${\cal S}_{\frac{5}{2}}(M;R) = R{\cal L}_{fr}/(I_{\frac{5}{2}})$
where $I_{\frac{5}{2}}$ is the submodule of $R{\cal L}_{fr}$
generated by the skein relation:\\
(i) \ \ $b_2L_2 + b_1L_1 + b_0L_0  + b_{\infty}L_{\infty} +
b_{-1}L_{-1} + b_{-\frac{1}{2}}L_{-\frac{1}{2}} = 0$,\\
 its mirror image:\\
($\bar{i}$) \ \ $b'_2L_2 + b'_1L_1 + b'_0L_0  + b'_{\infty}L_{\infty} +
b'_{-1}L_{-1} + b'_{-\frac{1}{2}}L_{-\frac{1}{2}} = 0$\\
 and the framing relation:\\
$L^{(1)} = a L$, where $a,b_2,b_2',b_{-\frac{1}{2}},
b'_{-\frac{1}{2}}$ are invertible elements in $R$ and 
$b_1,b'_1,b_0,b'_0$, $b_{-1},b'_{-1},b_{\infty}$, 
and $b'_{\infty}$ are any fixed elements of $R$. 
The links $L_2,L_1,L_0,L_{\infty},L_{-1},$ $L_{\frac{1}{2}}$ 
and $L_{-\frac{1}{2}}$
are illustrated in Fig.3.4.\footnote{Our notation is based on
Conway's notation for rational tangles. However it differs
from it by a sign. The reason is that the Conway
convention for a positive crossing is generally not
used in the setting of skein relations.}
\end{definition}
\centerline{\psfig{figure=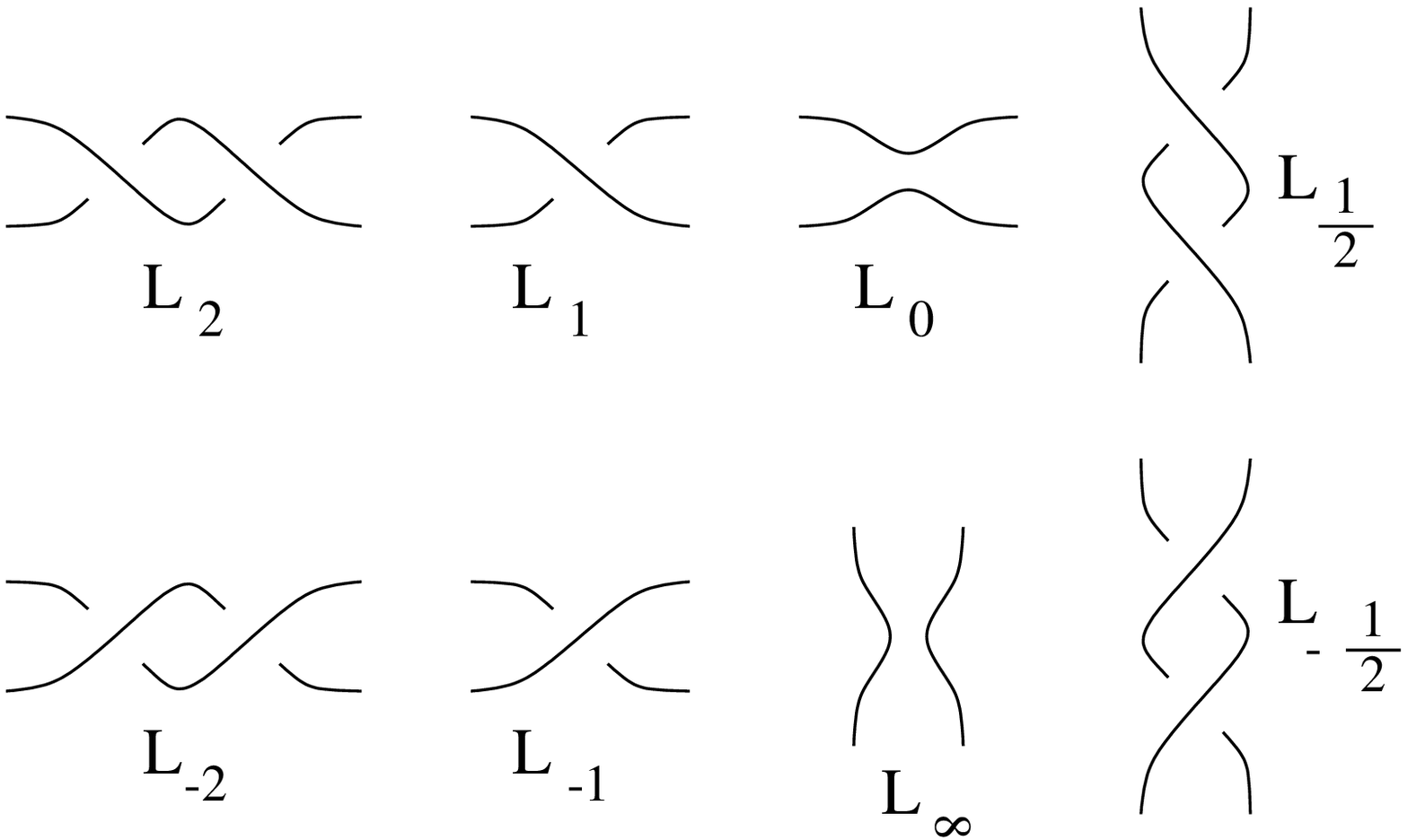,height=5.0cm}}
\centerline{Fig. 3.4}
If we rotate the figure from the relation (i) we obtain:\\
(i')  $b_{-\frac{1}{2}}L_2 + b_{-1}L_1 + b_{\infty}L_0  + b_{0}L_{\infty} +
b_{1}L_{-1} + b_{2}L_{-\frac{1}{2}} = 0$\\
One can use (i) and (i') to eliminate $L_{-\frac{1}{2}}$ and to get the relation:\\
$(b_2^2- b^2_{-\frac{1}{2}}) L_2 +  
(b_1b_2 - b_{-1}b_{-\frac{1}{2}})L_1 +
((b_0b_2 - b_{\infty}b_{-\frac{1}{2}})L_0 +
(b_{-1}b_2 - b_{1}b_{-\frac{1}{2}})L_{-1} +
(b_{\infty}b_2 - b_{0}b_{-\frac{1}{2}})L_{\infty}=0$.\\
Thus either we deal with the shorter relation (essentially that
of the fourth skein module described before) or all
coefficients are equal to 0 and therefore (assuming 
that there are no zero divisors in $R$)
$b_2 = \varepsilon b_{-\frac{1}{2}}$, $b_1 = \varepsilon b_{-1}$, 
and $b_0 = \varepsilon b_{\infty}$. Similarly we would get:
$b'_2 = \varepsilon b'_{-\frac{1}{2}}$, $b'_1 = \varepsilon b'_{-1}$, 
and $b'_0 = \varepsilon b'_{\infty}$. Here $\varepsilon = \pm 1$.
Assume for simplicity that $\varepsilon =1$. Further relations
among coefficients follows from the computation of the 
Hopf link component using the amphicheirality of the
unoriented Hopf link. 
Namely, by comparing diagrams in Figure 3.5 and its mirror image we get:
$$L\#H= -b_2^{-1}(b_1(a + a^{-1}) + a^{-2}b_2 + b_0(1+T_1))L$$
$$L\#H= -{b'}_2^{-1}(b'_1(a + a^{-1}) + a^{2}b'_2 + b'_0(1+T_1))L.$$
Possibly the above equalities give the only other relation among
coefficients (in the case of $S^3$), but I would present below
the simpler question (assuming $a=1, b_x=b'x$ and writting
$t^n$ for $T_n$).\\
\ \\
\centerline{\psfig{figure=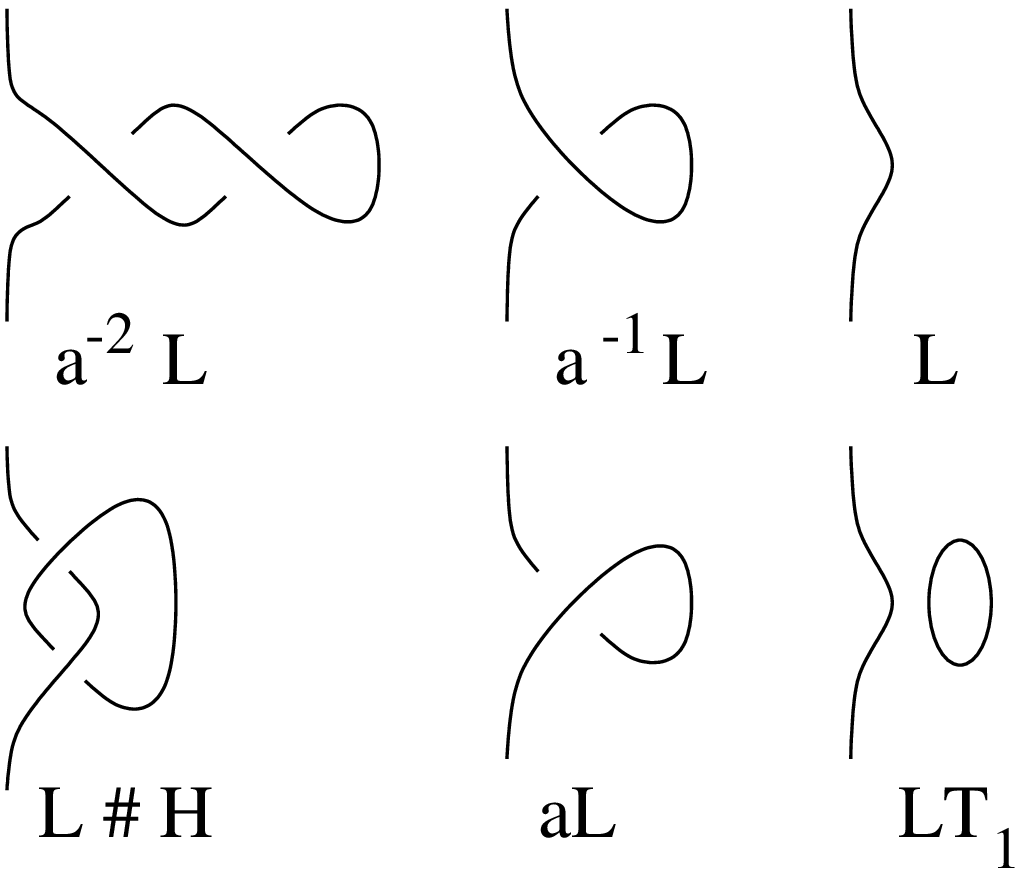,height=5.5cm}}
\centerline{Fig. 3.5}

\begin{question}
There is a polynomial invariant of unoriented links in $S^3$,
$P_{\frac{5}{2}}(L) \in Z[b_0,b_1,t]$ which satisfies:
\begin{enumerate}
\item [(i)] Initial conditions: $P_{\frac{5}{2}}(T_n) = t^n$, where $T_n$ is a trivial
link of $n$ components.
\item [(ii)] Skein relations
$$P_{\frac{5}{2}}(L_2) + b_1P_{\frac{5}{2}}(L_1) +  b_0P_{\frac{5}{2}}(L_0) + 
b_0P_{\frac{5}{2}}(L_{\infty}) + b_1P_{\frac{5}{2}}(L_{-1}) +  
P_{\frac{5}{2}}(L_{-\frac{1}{2}})
=0.$$
$$P_{\frac{5}{2}}(L_{-2}) + b_1P_{\frac{5}{2}}(L_{-1}) +  b_0P_{\frac{5}{2}}(L_0) + 
b_0P_{\frac{5}{2}}(L_{\infty}) + b_1P_{\frac{5}{2}}(L_{1}) +  
P_{\frac{5}{2}}(L_{\frac{1}{2}})
=0.$$ 
\end{enumerate}
\end{question}
Notice that by taking the difference of our skein relations one get
the interesting identity:
$$P_{\frac{5}{2}}(L_2)  - P_{\frac{5}{2}}(L_{-2}) =
P_{\frac{5}{2}}(L_{\frac{1}{2}}) - P_{\frac{5}{2}}(L_{-\frac{1}{2}}).$$
Nobody has yet seriously studied the skein module 
${\cal S}_{\frac{5}{2}}(M;R)$ so everything you can get will
be a new exploration, even a table of the polynomial $P_{\frac{5}{2}}(L)$
for small links, $L$.\\
I wish you luck!

\ \\
Department of Mathematics\\
George Washington University  \\
Washington, DC 20052 \\
USA\\
e-mail: przytyck@gwu.edu
\end{document}